\documentclass[1p]{elsarticle_modified}
\usepackage{amssymb,amsmath,amscd,tikz-cd,amsthm}
\usepackage{caption,subcaption}
\usepackage{url}
\usepackage[bookmarks=true,%
colorlinks=true,% 
linkcolor=blue,%
citecolor=blue,%
filecolor=blue,%
menucolor=blue,%
urlcolor=blue,%
breaklinks=true]{hyperref}
\usepackage{slashed}    
\usepackage{verbatim}
\usepackage{etoolbox}
\usepackage{mathtools}
\usepackage[normalem]{ulem}
\graphicspath{{figures/}}

\newtheorem{claim}{Claim}[section]

\usetikzlibrary{
	arrows,
	decorations.markings,
	positioning
}

\newcommand{\circled}[2][]{
	\tikz[baseline=(char.base)]{
		\node[shape = circle, draw, inner sep = 0.6pt]
		(char) {\phantom{\ifblank{#1}{#2}{#1}}};%
		\node at (char.center) {\makebox[0pt][c]{#2}};}}
\robustify{\circled}

\newenvironment{talign*}
{\let\displaystyle\textstyle\csname align*\endcsname}
{\endalign}

\def\sl{\mathrm{SL}_2(\mathbb{C})}
\def\psl{\mathrm{PSL}_2(\mathbb{C})}

\def\Hbb{\mathbb{H}}
\def\Cbb{\mathbb{C}}

\def\Ccal{\mathcal{C}}
\def\Bcal{\mathcal{B}}

\def\Tcal{\mathcal{T}}
\def\Dcal{\mathcal{D}}

\def\cl#1{{\small \circled[10]{#1}}}

\def\Hb{{\overline{\mathbb{H}^3}}}
\def\pmo{\{\pm 1\}}
\def\sl{\operatorname{\textup{SL}}(2,\Cbb)}
\def\psl{\operatorname{\textup{PSL}}(2,\Cbb)}
\def\pmo{\{\pm1\}}

\def\Hb{{\overline{\mathbb{H}^3}}}

\theoremstyle{definition}
\newtheorem{theorem}{Theorem}[section]
\newtheorem{proposition}[theorem]{Proposition}
\newtheorem{definition}[theorem]{Definition}
\newtheorem{example}[theorem]{Example}
\newtheorem{remark}[theorem]{Remark}

\begin{document}
	
	\begin{frontmatter}

		\title{Octahedral developing of knot complement II: \\  Ptolemy coordinates and applications}
		
		\author{Hyuk Kim\fnref{h_kim}}
		\fntext[hkim]{%	\address  
			Department of Mathematical Sciences, Seoul National University, Seoul, 08826, Korea.
			%\thanks
			Supported by Basic Science Research Program through the NRF of Korea funded by the Ministry of Education (NRF-2018R1A2B6005691).}
		\ead{hyukkim@snu.ac.kr}
		
		\author{Seonhwa Kim\fnref{skim}}
		\fntext[skim]{	%\address
			Department of Mathematics, 	University of Seoul, Seoul 02504, Korea.
			%\thanks
			Supported  by Basic Science Research Program through the
			National Research Foundation of Korea(NRF) funded by the Ministry of Education(2022R1I1A1A01063774) and
			by 
			the Institute for Basic Science (IBS-R003-D1).
		}
		\ead{seonhwa17kim@diagram.site}
		
		\author{Seokbeom Yoon\fnref{syoon}}
		\fntext[syoon]{	%\address
			International Center for Mathematics, Department of Mathematics,
			Southern University of Science and Technology, Shenzhen, China.
			%\thanks
		}
		\ead{sbyoon15@gmail.com}

		\begin{abstract}
			It is known that a knot complement (minus two points) decomposes into ideal octahedra with respect to a given knot diagram. 
			In this paper, we study the Ptolemy variety for such an octahedral decomposition in perspective of Thurston's gluing equation variety.
			More precisely, we compute explicit Ptolemy coordinates in terms of segment and region variables, the coordinates of the gluing equation variety motivated from the volume conjecture.
			As a consequence, we present an explicit formula for computing the obstruction to lifting a boundary-parabolic $\psl$-representation to boundary-unipotent $\sl$-representation.
			%	$(\mathrm{PSL}(2,\mathbb{C}),P)$-representation of the knot group to a $(\mathrm{SL}(2,\mathbb{C}),P)$-representation.	
			We also present a diagrammatic algorithm to compute a holonomy representation of the knot group.

		\end{abstract}
		
	\end{frontmatter}

	\tableofcontents

\section{Introduction}\label{sec:Intro}

\subsection{Gluing equation varieties and Ptolemy varieties}
Let $M$ be a compact orientable 3-manifold with non-empty boundary and $\Tcal$ be an ideal triangulation of its interior.
The gluing equations of $\Tcal$ with cusp conditions define an algebraic set $V(\Tcal)$, called the \emph{gluing equation variety}, whose coordinates are the shapes of ideal tetrahedra of $\Tcal$ \cite{thurston1979geometry, neumann1985volumes}.
In addition, a pseudo-developing map  associates each point of $V(\Tcal)$ to a holonomy representation $\rho : \pi_1(M)\rightarrow \psl$, which is boundary-parabolic, up to conjugation (see, for instance,  \cite[Section 2.5]{kim2018octahedral}). 
Namely, there is a map 
\begin{equation*}
	\mathrm{hol} : V(\Tcal)\rightarrow \left\{ \begin{array}{c}\textrm{boundary-paraoblic} \\  \rho : \pi_1(M)\rightarrow \psl  \end{array} \right\}/_\textrm{Conjugation}.
\end{equation*}
On the other hand, inspired by \cite{fock2006moduli, zickert2009volume}, Garoufalidis-Thurston-Zickert \cite{garoufalidis2015complex} introduced another construction of algebraic sets called Ptolemy varieties.
The \emph{Ptolemy variety} $P^\sigma(\Tcal)$ of $\Tcal$ is defined for each class $\sigma\in H^2(M,\partial M;\pmo)$ and admits a  map
\begin{equation*} 
	\mathrm{hol} : P^\sigma(\Tcal)\rightarrow \left\{ \begin{array}{c}  \textrm{boundary-paraoblic} \\  \rho : \pi_1(M)\rightarrow \psl \\ \textrm{whose obstruction class is } \sigma \end{array} \right\}/_\textrm{Conjugation}
\end{equation*}
where the obstruction class of $\rho$ means the obstruction to lifting $\rho$ to a boundary-unipotent $\sl$-representation.

The way in which the Ptolemy variety produces holonomy representations is more practical, compared to the gluing equation variety, in the sense that there is no need to construct pseudo-developing maps.
Moreover, one can say that the Ptolemy variety is more \emph{efficient} than the gluing equation variety, as its coordinates directly produce the shapes of the ideal tetrahedra. Indeed, there is a surjective regular map $\psi$ defined tetrahedron-wise with the following diagram (\cite[Theorem~1.12]{garoufalidis2015complex}):
\begin{equation} \label{eqn:psi}
	\begin{tikzcd}[column sep=20pt]
		\displaystyle P(\Tcal) := \coprod_{\sigma}  P^\sigma(\Tcal) \arrow[rd,"\mathrm{hol}"] \arrow[d, two heads, swap, "\psi"] & \\
		V(\Tcal) \arrow[r, "\mathrm{hol}"] &  \left\{ \begin{array}{c}\textrm{boundary-paraoblic} \\  \rho : \pi_1(M)\rightarrow \psl  \end{array} \right\}/_\textrm{Conjugation}.			
	\end{tikzcd}
\end{equation}
%Remark that in general, it is complicated to compute the coordinates of the Ptolemy variety from the shapes of ideal tetrahedra, as the holonomy action should be involved in the computation.

\subsection{Main results}

For a knot $K \subset S^3$ it is known that the knot complement  with two points removed decomposes into ideal octahedra canonically with respect to a given knot diagram. Hence we obtain an ideal triangulation $\Tcal$ of $S^3 \setminus (K \cup \{p,q\})$, where $p \neq q \in S^3$ are two points not in $K$, by dividing each ideal octahedron to ideal tetrahedra.
In \cite{kim2018octahedral}, we have studied basic properties of such ideal triangulations, mainly focused on pseudo-developing maps and holonomy representations via the gluing equations.  
In particular, we have expressed the gluing equation variety of $\Tcal_4$ and $\Tcal_5$ in terms of so-called $z$- and $w$-variables, which were derived from the optimistic limit of the Kashaev invariant and colored Jones polynomials (see \cite{thurston1999hyperbolic, yokota_volume_2000, cho2013optimistic,cho2014optimistic} for details, and see also \cite{kim2018octahedral} for geometric approach). 
These variables defined originally as the formal parameters in the integral of q-series  also appear recently in a rather different quantum setting \cite{mcphail-snyder_hyperbolic_2022}. It is quite surprising to notice that these quantum variables can be interpreted in the hyperbolic geometry terms directly, which strongly suggests rather direct connection between these two areas. 

As a sequel to \cite{kim2018octahedral}, it is shown that such variables are practically helpful in connecting  geometry of knot complements and combinatorics of knot diagrams from the perspective of Ptolemy coordinates.
Precisely, 
a main object of this paper is the Ptolemy variety for the above ideal triangulation and  our goal is to construct an explicit map $\phi:V(\Tcal)\rightarrow P(\Tcal)$, a section of the surjective map $\psi : P(\Tcal) \rightarrow V(\Tcal)$ in the diagram~\eqref{eqn:psi} for the case of $\Tcal = \Tcal_4$ or $\Tcal_5$: 
\begin{equation} \label{eqn:phi}
	\begin{tikzcd}[column sep=20pt]
		\displaystyle P(\Tcal)  \arrow[rd,"\mathrm{hol}"] \arrow[d,  shift right=1, two heads, swap, "\psi"] & \\
		V(\Tcal) \arrow[u, shift right=1, swap, "\phi"] \arrow[r, "\mathrm{hol}"] &  \left\{ \begin{array}{c}\textrm{boundary-paraoblic} \\  \rho : \pi_1(M)\rightarrow \psl  \end{array} \right\}/_\textrm{Conjugation}.			
	\end{tikzcd}
\end{equation}

In general, the construction of a section $\phi : V(\Tcal) \rightarrow P(\Tcal)$ is well-known \cite{zickert2009volume}, but it can not be defined tetrahedron-wise (contrast to $\psi$) as the construction of a pseudo-developing map involves. 
However, surprisingly, it can be written ``almost octahedron-wise'' in our case: Ptolemy coordinates near a crossing are given by $z$- or $w$-variables around the crossing with a single scaling parameter.
This fact will be fruitfully used for our formulas and their proofs.

As a result, we obtain the Ptolemy coordinates (the coordinates of the Ptolemy variety) explicitly  in terms of $z$- or $w$-variables.
As we mentioned earlier, the Ptolemy coordinates efficiently produce holonomy representations, hence the map $\phi$ allows us to compute boundary-parabolic representations $\rho$ of the knot group directly  in terms of $z$- or $w$-variables. 
This leads us to several diagrammatic formulas related to $\rho$. As examples, in this paper, we present diagrammatic formulas for computing  

\begin{itemize}
		\item the obstruction class of $\rho$ (Theorem \ref{thm:obs});
		\item the cusp shape of $\rho$ (Theorem \ref{thm:longitude});
		\item the $\rho$-image of Wirtinger generators (Theorem \ref{thm:w_wirtinger}).
\end{itemize}

The second formula is computing 
the cusp shape of boundary-parabolic representation which is defined as the longitude translation divided by the meridian translation, and  determines the Euclidean structure of the cusp cross-section if $\rho$ is the holonomy of a complete hyperbolic metric.  In the study of  the Volume conjecture,  another cusp shape formula was also obtained  by Y. Yokota~\cite{yokota_cusp_2016}, that is given in terms of  the Hessian of volume potential function. 
For now, it seems unclear whether our cusp shape formula can be obtained in a similar way. 

We organize the paper as follows.
In Section~\ref{sec:deco}, we recall  how Ptolemy coordinates are computed from a pseudo-developing map in general. 
In Section~\ref{sec:main}, we present explicit Ptolemy coordinates for our case in terms of $z$- and $w$-variables. This will be done in two steps: we first compute Ptolemy coordinates for each ideal octahedron and then glue them compatibly.
In Sections~\ref{sec:obs} and \ref{sec:wirtgen}, we prove the diagrammatic formulas that we mentioned above.

For simplicity of exposition, we restrict our attention to knots in $S^3$. However, most of the discussions in this paper also work for links in $S^3$. 
For example, 
the obstuction class and the longitude for a link are established by multiplying $\sigma(c_i)$ (Theorem \ref{thm:obs}) 
and adding $\lambda(c_i)$ (Theorem  
\ref{thm:longitude})   whenever the crossing is under-passed along the corresponding link-component, since they are defined only link-component-wisely.

\section{Ptolemy coordinates from developing maps}\label{sec:deco}

We devote this section to recall some notions and facts in \cite{zickert2009volume, garoufalidis2015complex, garoufalidis2015ptolemy}.
We reorganize them from the point of view of a developing map so that they fit more naturally to this paper.

\subsection{Developing maps and decorations}

Let $M$ be a compact orientable 3-manifold with non-empty boundary and
$G$ be either $\psl$ or $\sl$. 
Throughout the section, we fix a $(G,P)$-representation $\rho : \pi_1(M)\rightarrow G$ defined as follows.
\begin{definition}
	A representation $\rho:\pi_1(M) \rightarrow G$ is called a \emph{$(G,P)$-representation} if up to conjugation $\rho(\pi_1(\Sigma))$ lies in $P$  for each boundary component $\Sigma$ of $M$. Here $P$ is the subgroup consisting of upper triangular matrices with ones on the diagonal.   
	A $(G,P$)-representation is also called \emph{boundary-paraoblic} if $G=\psl$ and \emph{boundary-unipotent} if $G=\sl$.  
\end{definition}

We denote by $N$ the universal cover of $M$, $\hat{N}$ the space obtained from $N$ by collapsing each boundary component of $N$ to a point, and  $v(\hat{N})$ the set of collapsed points.
Note that the deck transformation action of $\pi_1(M)$ on $N$ induces the $\pi_1(M)$-action on $\hat{N}$ and $v(\hat{N})$.
\begin{definition}\label{def:bdry_para}
	A \emph{(pseudo-)developing map} $\mathfrak{D}:\hat{N} \rightarrow \overline{\Hbb^3}$ is a $\rho$-equivariant map such that $\mathfrak{D}(x) \in \partial \Hb$ if and only if $x \in v(\hat{N})$. Here $\rho$-equivariance means that $\mathfrak{D}(\gamma \cdot x) = \rho(\gamma)  \mathfrak{D}(x)$ for $\gamma \in \pi_1(M)$ and $x \in \hat{N}$.
\end{definition} 

\begin{definition}
	A \textit{decoration} $\Dcal : v(\hat{N}) \rightarrow G/P$  is a $\rho$-equivariant map where $G/P$ is the left $P$-coset space.
\end{definition}

The above two notions are related as follows.
For a $P$-coset $gP$ we define 
\[gP(\infty):= g'(\infty) \in \partial \Hb=\Cbb \cup \{\infty\}\]
for any coset representative $g'$ of $gP$, viewed as a M\"obious transformation on $\partial \overline{\Hbb^3}$. Note that $gP(\infty)$ is well-defined, as every element of $P$ fixes $\infty$.

\begin{proposition} \label{prop:dev_dec}
	For a developing map $\mathfrak{D}:\hat{N} \rightarrow \overline{\Hbb^3}$ there is a decoration $\Dcal  : v(\hat{N}) \rightarrow G/P$ such that 
	\begin{equation} \label{eqn.prop}
		(\Dcal(x))(\infty)=\mathfrak{D}(x) \quad \textrm{for all } x \in v(\hat{N}).
	\end{equation}
	Namely, the following diagram commutes:
	\begin{equation*}
		\begin{tikzcd}[column sep=30pt]
			& G/P \arrow[d]\\
			v(\hat{N})\arrow[ru,"\Dcal"] \arrow[r, swap,"\mathfrak{D}|_{v(\hat{N})}"] & \partial\Hb			
		\end{tikzcd}
	\end{equation*}
	where the vertical map is given by evaluating $\infty$.
\end{proposition}
\begin{proof} 
	Let $V=\{v_1,\ldots,v_h\} \subset v(\hat{N})$ be the set of  representatives of the $\pi_1(M)$-orbits in $v(\hat{N})$ (so $h$ is the number of boundary components of $M$).  
	For each $v_i\in V$ we define $\Dcal(v_i):=g_iP$ where $g_i \in G$ is any element satisfying $g_i(\infty)=\mathfrak{D} (v_i)$. 
	We then define $\Dcal(v)$ for other $v\in v(\hat{N})$ not in $V$, $\rho$-equivariantly: \[ \Dcal(v):=\rho(\gamma) \Dcal(v_i)= \rho(\gamma)g_iP\]
	for any pair of $v_i\in V$ and $\gamma \in \pi_1(M)$ such that $v=\gamma \cdot v_i$.
	Such  a pair clearly exists, but may not be unique: choosing $v_i\in V$ is unique but not for $\gamma \in \pi_1(M)$.  
	We thus need to check well-definedness. Suppose that $\gamma \cdot v_i = \gamma' \cdot v_i$ for some  $\gamma$ and $\gamma' \in \pi_1(M)$. Then $\rho(\gamma^{-1}\gamma')$ is a parabolic element fixing $g_i(\infty)=\mathfrak{D}(v_i)$ and  thus an element of $g_i P g_i^{-1}$. Note that $g_i P g_i^{-1}$ is the set of all parabolic elements of $G$ fixing $\mathfrak{D}(v_i)$. It follows that $\rho(\gamma) g_iP=\rho(\gamma') g_iP$. This proves that $\Dcal : v(\hat{N}) \rightarrow G/P$ is well-defined. 
	The fact that $\Dcal$ satisfies the equation~\eqref{eqn.prop} is clear from its construction.
\end{proof}

The notion of a decoration admits a nice geometric interpretation (see \cite[\S 3.1]{zickert2009volume}), but it would not be considered in this paper. 
%Also, the construction of $\Dtilde$ in Proposition~\ref{prop:dev_dec} is not unqiue; there are $(\Cbb^\times)^h$ many choices. This implies that the map $\psi$ in the diagram \eqref{eqn:psi} is $(\Cbb^\times)^h$-to-one  (see Theorem 1.12 of \cite{garoufalidis2015complex}).

\subsection{Ptolemy varieties} \label{sec:ptol}

In this section, we recall the definition of Ptolemy varieties and explain how their coordinates are computed from a developing map.

Let $\Tcal$ be an ideal triangulation of the interior of $M$ with edge set $e(\Tcal)$ and face set $f(\Tcal)$.
We denote by $Z^2(\Tcal;\pmo)$ the set of all maps $\sigma : f(\Tcal) \rightarrow \pmo$ such that \[\sigma(f_0)\sigma(f_1)\sigma(f_2)\sigma(f_3) =1\] if $f_0,\ldots,f_3$ are the faces of an ideal tetrahedron of $\Tcal$. An element of $Z^2(\Tcal;\pmo)$ can be viewed as a usual 2-cocycle, an element of $Z^2(M,\partial M;\pmo)$, hence defines a class in $H^2(M,\partial M;\pmo)$.

\begin{definition} \label{def:ptolemy}
	A \emph{Ptolemy assignment with the obstruction cocycle $\sigma \in Z^2(\Tcal;\{\pm1\})$} is a set map $c:e(\Tcal)\rightarrow \Cbb^\ast$  such that
	$-c(e)=c(-e)$ for all $e \in e(\Tcal)$ and
	\begin{equation} \label{eqn:ptolemy_obs}
		\sigma(f_2)c(l_{02})c(l_{13})= \sigma(f_3)c(l_{03})c(l_{12})+ \sigma(f_1) c(l_{01}) c(l_{23})
	\end{equation}
	for each ideal tetrahedron $\Delta$ of $\Tcal$. Here  $-e$ denotes the same edge $e$ with the opposite orientation, $l_{ij}$ is the oriented edge of $\Delta= [v_0,v_1,v_2,v_3]$ running from $v_i$ to $v_j$, and $f_i$ is the face of $\Delta$  opposite to $v_i$ (see Figure \ref{fig:trunctaion}).
	The \emph{Ptolemy variety} $P^\sigma(\Tcal)$ is the set of all Ptolemy assignments with the obstruction cocycle $\sigma$.
	\begin{figure}[!h]
		\centering
		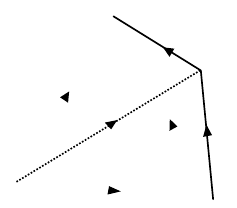
		\caption{An ideal tetrahedron.}
		\label{fig:trunctaion}
	\end{figure}
\end{definition}

The ideal triangulation $\Tcal$ endows $M$ with a decomposition into truncated tetrahedra whose triangular faces triangulate the boundary $\partial M$.  We denote by $v(M)$ and $e(M)$ the vertex set and the edge set of the decomposition of $M$, respectively. We call an edge of $M$ a \emph{short-edge} if it is contained in $\partial M$, and a \emph{long-edge}, otherwise.
This decomposition lifts to the universal cover $N$ of $M$, hence the vertex set $v(N)$ and the edge set $e(N)$ (as well as the notion of short- and long-edges) of $N$ are defined.

%%%%%%Let $\mathfrak{D}:\hat{N} \rightarrow \overline{\Hbb^3}$  be a developing map with a decoration $\Dcal : v(\hat{N}) \rightarrow G/P$ described as in Proposition~\ref{prop:dev_dec}.
%%%%%%%We assume that every ideal tetrahedron of $\Tcal$ has a non-degenerate image under the developing map, i.e., $\mathfrak{D}(v_0),\ldots,\mathfrak{D}(v_3)$ are mutually distinct if $v_0,\ldots,v_3$ are the vertices of an ideal tetrahedron.
Note that if we start from a solution of the gluing equations, then its associated developing map automatically satisfies the above assumption.	

It is proved in \cite[Lemma 3.5]{zickert2009volume} that there is a unique $\rho$-equivariant map $\Ccal : v(N) \rightarrow G$ such that 
\begin{itemize}
	\item \label{eqn:con1} $\Ccal(x)  \in \Dcal(v)$ if $x$ is in the boundary component of $N$ collapsed to $v\in v(\hat{N})$;
	\item \label{eqn:con2} $\Ccal(x_0)^{-1} \Ccal(x_1)$ is a counter-diagonal matrix  if $x_0$ and $x_1$ are joined by a long-edge of $N$.
\end{itemize}  
We define a map $\Bcal : e(N) \rightarrow G$ by letting
\begin{equation} \label{eqn:bcal}
	\Bcal(e):=\Ccal (x_0)^{-1} \Ccal (x_1)
\end{equation} for $e=[x_0,x_1] \in e(N)$ where $x_0$ and $x_1$ are the initial and terminal vertices of $e$, respectively.
Then for $\gamma \in \pi_1(M)$ and $e =[x_0,x_1] \in e(N)$ we have
\begin{align*}
	\Bcal(\gamma \cdot e) &= \Ccal(\gamma \cdot x_0)^{-1} \Ccal(\gamma \cdot x_1)\\ 
	&= (\rho(\gamma) \Ccal(x_0))^{-1} \rho(\gamma) \Ccal(x_1)\\
	&= \Ccal(x_0)^{-1} C(x_1) \\
	&= \Bcal(e).
\end{align*} 
Hence the map $\Bcal$ projects down to  $\Bcal:e(M) \rightarrow G$ (also denoted by $\Bcal$).
The definition~\eqref{eqn:bcal} implies that $\Bcal$ automatically satisfies the cocycle condition:
\begin{itemize}
	\item  $\Bcal(e) \Bcal(-e)=I$ for $e \in e(M)$;
	\item  $\Bcal(e_1)\cdots\Bcal(e_m)=I$ if $e_1,\ldots,e_m$ are the boundary edges of a face of $M$  in order ($m$ is either 3 or 6, as a face of $M$ is either triangular or hexagonal).
\end{itemize}
We obtain a representation $\rho_0: \pi_1(M)\rightarrow G$ induced from the cocycle $\Bcal$. Precisely, for $\gamma \in\pi_1(M)$  we homotope it to an edge-path in $M$ and define $\rho_0(\gamma)$ by the product of the matrices assigned by $\Bcal$ along the edge-path. 
Note that $\rho_0$ is well-defined, as $\Bcal$ satisfies the cocycle conditions.

\begin{proposition} \label{prop:hol_recv}
	The representation $\rho_0$ is conjugate to $\rho$, the one that we started with.
\end{proposition}        
\begin{proof} 
	We choose a base point $x$ of $\pi_1(M)$  in $v(M)$ together with its lifting $\widetilde{x}\in v(N)$. For any loop $\gamma \in \pi_1(M,x)$ let $\gamma_0$ be an edge-path in $M$ based at $x$ and homotopic to $\gamma$. Recall that $\rho_0(\gamma)$ is given by the product of the matrices assigned by $\Bcal$ along $\gamma_0$, which is clearly equal to the product of the matrices along the lifting $\widetilde{\gamma}_0$ of $\gamma_0$. Here $\widetilde{\gamma}_0$ is an edge-path in $N$ based at $\widetilde{x}$. As  $\widetilde{\gamma}_0$ runs from $\widetilde{x}$ to $\gamma \cdot \widetilde{x}=\gamma_0 \cdot \widetilde{x}$, we have
	\begin{equation}
		\rho_0(\gamma)=\Ccal(\widetilde{x})^{-1} \Ccal(\gamma \cdot \widetilde{x})=\Ccal(\widetilde{x})^{-1}\rho(\gamma) \Ccal(\widetilde{x})
	\end{equation}
	which proves that $\rho_0$ and $\rho$ are conjugate.
\end{proof}

The cocycle $\Bcal : e(M) \rightarrow G$ is called a \emph{natural $(G,P)$-cocycle} (or simply a \emph{natural coycle}) in \cite{garoufalidis2015gluing, garoufalidis2015ptolemy},  since it takes matrices of particular forms: $\Bcal(e) \in P$ if $e$ is a short-edge and  $\Bcal(e)$ is counter-diagonal if $e$ is a long-edge.

When $G=\sl$, a natural cocycle is reduced to a map $c : e(M) \rightarrow \Cbb$ by letting  \begin{equation}\label{eqn:assigning}
	\left \{
	\begin{array}{l}
		\Bcal(e)=\arraycolsep=1pt\begin{pmatrix} 0 & -c(e)^{-1} \\ c(e) & 0 \end{pmatrix}  \textrm{ for all long-edges } e,  \\[15pt]
		\Bcal(e)=\begin{pmatrix} 1 & c(e) \\ 0 & 1  \end{pmatrix} \textrm{ for all short-edges } e.
	\end{array} \right.
\end{equation}
According to $e\in e(M)$, we call $c(e)$ either a \textit{short-edge parameter} or a \textit{long-edge parameter}.
It is proved in \cite[Lemma 3.3]{zickert2009volume} that the long-edge parameters determine all short-edge parameters.
Forgetting all short-edge parameters and identifying each long-edge of $M$ with an edge of $\Tcal$ in an obvious way, we obtain a map (also denoted by $c$) \[c : e(\Tcal) \rightarrow \Cbb^\ast.\] 
Note that a long-edge parameter can not be zero (see the equation~\eqref{eqn:assigning}).
It is proved in \cite[Lemma 3.5]{zickert2009volume} that 
the map $c$ can be directly computed from the decoration $\Dcal$ as follows.
\begin{equation}\label{eqn:detf}
	c(e) = \det \left(\Dcal(v_0) \binom{1}{0},\, \Dcal(v_1) \binom{1}{0}\right), \quad e \in e(\Tcal)
\end{equation}
where ${v}_0$  and ${v}_1$ are the initial and terminal vertices of any lifting of $e$, respectively.
Note that the above determinant does not depend on the choice of the lifting, as $\Dcal$ is $\rho$-equivariant, and that for a $P$-coset $gP$, $gP\binom{1}{0}$ means a vector $g'\binom{1}{0} \in \Cbb^2$ where $g'$ is any coset representative of $gP$.

\begin{proposition} \label{prop:ptol}
	We have $c \in P^\sigma(\Tcal)$ for the trivial obstruction cocycle $\sigma$. Here we say that $\sigma \in Z^2(\Tcal;\pmo)$ is trivial if $\sigma(f) =1$ for all $f \in f(\Tcal)$.
\end{proposition}
\begin{proof}
	For an ideal tetrahedron $\Delta$ of $\Tcal$ with vertices $v_0,\ldots,v_3$ let $V_i:=\Dcal(v_i) \binom{1}{0}$ for $0 \leq i \leq 3$. 
	Plugging the equation \eqref{eqn:detf} to the equation \eqref{eqn:ptolemy_obs} with the trivial cocycle $\sigma$, we obtain  
	\[\textrm{det}(V_0,V_2) \textrm{det}(V_1,V_3) = \textrm{det}(V_0,V_3)\textrm{det}(V_1,V_2) + \textrm{det}(V_0,V_1) \textrm{det}(V_2,V_3),\]
	which is true due to the Pl\"{u}cker relation. 
	Also,  the equation \eqref{eqn:detf} implies that $-c(e)=c(-e)$ for $e \in e(\Tcal)$, hence $c \in P^\sigma(\Tcal)$ for the trivial cocycle $\sigma$.
\end{proof}

For $G=\psl$ we similarly define short-edge parameters and long-edge parameters: the former is well-defined, but the latter has sign-ambiguity, as the counter-diagonal matrix in the equation \eqref{eqn:assigning} is the same $\psl$-matrix for $c(e)$ and $-c(e)$. 
We thus need to choose a sign for each long-edge parameter in order  to define a map $c:e(\Tcal)\rightarrow \Cbb^\ast$. 
%Such a set map $c$ may not be a point of $P^\sigma(\Tcal)$ for the trivial $\sigma$ but is a point of $P^\sigma(\Tcal)$ for some $\sigma \in Z^2(\Tcal;\pmo)$.
% 	(see Proposition \ref{prop:ptolob} below). 
%	To give a precise statement, we first recall the notion of obstruction classes for a $\psl$-representation.

% 	A representation $\pi_1(\Mline)\to\textrm{PSL}(2,\mathbb{C})$ may or may not lift to $\textrm{SL}(2,\mathbb{C})$ and the obstruction to such a lifting is a class in $H^2(\Mline;\{\pm 1\})$.
% 	Similarly, by the \emph{obstruction class of a $(\psl,P)$-representation $\pi_1(\Mline)\rightarrow \psl$} we mean the obstruction to lifting the representation to a $(\sl,P)$-representation. This is a class in $H^2(\Mline,\partial \Mline;\pmo) \simeq H^2(\Mhat;\pmo)$. 
% 	We refer to \cite{garoufalidis2015complex}, \cite[\S2.1]{CYZ2019} for details (see also the proof of Proposition \ref{prop:ptolob}).

\begin{proposition} \label{prop:ptolob}
	We have $c \in P^\sigma(\Tcal)$ for some $\sigma \in Z^2(\Tcal;\{\pm1\})$ whose class in $H^2(M,\partial M;\{\pm1\})$ is the obstruction class of $\rho$.
\end{proposition} 
\begin{proof} 
	Once the map $c:e(\Tcal)\rightarrow \Cbb^\ast$ is given, we have a unique $\sl$-lift $\widetilde{\Bcal}: e(M) \rightarrow \sl$ of the natural $(\psl,P)$-cocycle $\Bcal$ such that $\widetilde{\Bcal}(e) \in P$  if $e$ is a short-edge and 
	\begin{equation}
		\widetilde{\Bcal}(e) = \begin{pmatrix} 0 & -c(e)^{-1} \\ c(e) & 0 \end{pmatrix} 
	\end{equation}
	if $e$ is a long-edge of $M$. Clearly, $\widetilde{\Bcal}$ may not be a cocycle. 
	The product of $\widetilde{\Bcal}$-matrices along the boundary of a triangular face of $M$ is $I$, but this product along the boundary of a hexagonal face is  either $I$ or $-I$. 
	This defines a 2-cocycle \[\sigma \in Z^2(M,\partial M;\pmo) \simeq Z^2(\Tcal;\pmo)\] whose class in $H^2(M,\partial M;\pmo)$ is by definition the obstruction class of $\rho$. 
	We refer to \cite[\S2]{CYZ2019} for details.
	
	We now consider an ideal tetrahedron $\Delta=[v_0,v_1,v_2,v_3]$ of $\Tcal$ with its truncation $\overline{\Delta}$. We denote by $l_{ij}$ the long-edge of $\overline{\Delta}$ running from $v_i$ to $v_j$ and $s^{i}_{jk}$ the short-edge of $\overline{\Delta}$ running from $l_{ij}$ to $l_{ik}$. For simplicity we let $\widetilde{\Bcal}_{ij}:=\widetilde{\Bcal}(l_{ij})$, $\widetilde{\Bcal}^i_{jk}:=\widetilde{\Bcal}(s^i_{jk})$, and 
	\[V_0 := \binom{1}{0}, \quad V_1:=\widetilde{\Bcal}_{01} \binom{1}{0}, \quad V_i:= \widetilde{\Bcal}^0_{1i} \widetilde{\Bcal}_{0i}\binom{1}{0}, \quad i=2,3.\]
	Since every $\widetilde{\Bcal}^i_{jk} \in P$ fixes $V_0$, we have $\det (V_0,V_i)=c(l_{0i})$ for $1 \leq i \leq 3$. 
	Also, we have
	\begin{align*}
		\det(V_1,V_2)&= \textrm{det} \left(\widetilde{\Bcal}_{01} \binom{1}{0},\, \widetilde{\Bcal}^0_{12}\widetilde{\Bcal}_{02} \binom{1}{0}\right) \\
		&= \det \left(\widetilde{\Bcal}_{01} \binom{1}{0},\, \sigma(f_3) \widetilde{\Bcal}_{01}\widetilde{\Bcal}^1_{02}\widetilde{\Bcal}_{12}\widetilde{\Bcal}^2_{10} \binom{1}{0}\right) \\
		&= \sigma(f_3) \det\left(\binom{1}{0},\, \widetilde{\Bcal}^1_{02}\,\widetilde{\Bcal}_{12}\binom{1}{0} \right)\\
		&=\sigma(f_3)  c(l_{12}).
	\end{align*}
	Here $f_i$ is the hexagonal face of $\overline{\Delta}$  opposite to $v_i$. Note that we use the definition of $\sigma$ for the second equality, that is, \[\widetilde{\Bcal}_{01} \mkern 2mu \widetilde{\Bcal}^1_{02} \mkern 2mu  \widetilde{\Bcal}_{12} \mkern 2mu  \widetilde{\Bcal}^2_{10} \mkern 2mu  \widetilde{\Bcal}_{20} \mkern 2mu  \widetilde{\Bcal}^0_{21}= \sigma(f_3) I.\]
	A similar computation shows that the Pl\"{u}cker relation for $V_0, \ldots,V_3$ is equivalent to the equation \eqref{eqn:ptolemy_obs}:
	\begin{align*}
		&\det(V_0,V_2) \det(V_1,V_3) = \det(V_0,V_3)\det(V_1,V_2) + \det(V_0,V_1) \det(V_2,V_3) \\
		\Leftrightarrow\quad & \sigma(f_2) c(l_{02})  c(l_{13}) = \sigma(f_3) c(l_{03})  c(l_{12}) + \sigma(f_1)  c(l_{01}) c(l_{23}).
	\end{align*}
	This proves that $c \in P^\sigma(\Tcal)$.
\end{proof}

\begin{remark}
	If we change the sign-choices that we made to define the map $c:e(\Tcal)\rightarrow \Cbb^\ast$, then the obstruction cocycle $\sigma$ in Proposition~\ref{prop:ptolob} would change, but its class in $H^2(M,\partial M;\pmo)$ should be the same. Moreover, it is known that 
	two Ptolemy varieties $P^{\sigma_0}(\Tcal)$ and $P^{\sigma_1}(\Tcal)$ are canonically isomorphic if  obstruction cocycles $\sigma_0$ and $\sigma_1$ represent the same class in $H^2(M,\partial M;\pmo)$.
	Thus the way of choosing the signs is not important.		
\end{remark}

\section{Ptolemy coordinates for an octahedral decomposition}\label{sec:main}

\begin{figure}[h!]
	\centering
	\tikzset{
		commutative diagrams/.cd,
		arrow style=tikz,
		diagrams={>=stealth}}
	\begin{tikzcd}[column sep=2em]	
		\vcenter{\hbox{\includegraphics{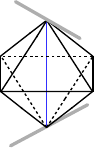}}} \ar[r]  &  \vcenter{\hbox{\includegraphics{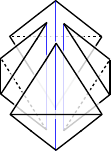}}} &  & 
		\vcenter{\hbox{\includegraphics{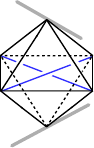}}}\ar[r]  &  \vcenter{\hbox{\includegraphics{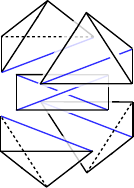}}}
	\end{tikzcd}
	\caption{Subdivisions of an ideal octahedron.}
	\label{fig:four_five}
\end{figure}

Let $K$ be a knot in $S^3$ with a diagram $D$. The complement $S^3 \setminus (K\cup\{p,q\})$ decomposes into ideal octahedra (one per a crossing) where $p\neq q\in S^3$ are two points not in $K$. See, e.g., \cite{thurston1999hyperbolic, weeks2005computation, kim2018octahedral}. We denote these ideal octahedra placed at the crossings by $o_1,\ldots,o_n$ (so $n$ is the number of crossings of $D$).
Subdividing each octahedron into ideal tetrahedra as in Figure \ref{fig:four_five}~(left) (resp., (right)), we obtain an ideal triangulation $\Tcal_4$ (resp., $\Tcal_5$) of $S^3 \setminus (K\cup\{p,q\})$. Figure \ref{fig:four_five}  is taken from \cite{kim2018octahedral}.

\begin{figure}[!h] 
	\centering  	
	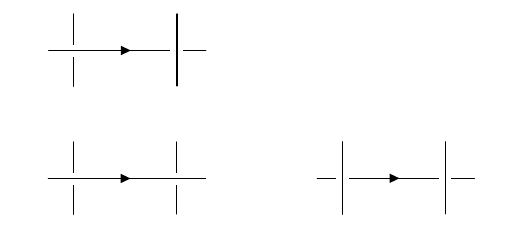
	\caption{Segment variables around a segment.}
	\label{fig:local_diagram_z}
\end{figure}

In \cite{kim2018octahedral} we expressed the gluing equations of $\Tcal_4$ in terms of \emph{$z$-variables}, also called \emph{segment variables}. As the name suggests, they are non-zero variables assigned to segments of $D$, where a segment means an edge of a diagram, viewed as a 4-valent graph. Precisely, the gluing equations of $\Tcal_4$ are given as:
\begin{equation} \label{eqn:hyp_z}
	\left\{
	\begin{array}{cccc}
		\dfrac{z_c-z_a}{z_c-z_b} \cdot \dfrac{z_d(z_c-z_e)}{z_e(z_c-z_d)} &=& 1 & \textrm{for Figure \ref{fig:local_diagram_z}~(a)} \\[10pt]
		\dfrac{z_a(z_c-z_b)}{z_b(z_c-z_a)} \cdot \dfrac{z_c-z_d}{z_c-z_e}&=&1 & \textrm{for Figure \ref{fig:local_diagram_z}~(b)} \\[10pt]
		\dfrac{z_c-z_a}{z_c-z_b} \cdot \dfrac{z_c-z_e}{z_c-z_d} &=& 1 & \textrm{for Figure \ref{fig:local_diagram_z}~(c)} \\[10pt]
		\dfrac{z_a(z_c-z_b)}{z_b(z_c-z_a)} \cdot \dfrac{z_e(z_c-z_d)}{z_d(z_c-z_e)}    &=& 1 & \textrm{for Figure \ref{fig:local_diagram_z}~(d)}
	\end{array}
	\right.
\end{equation}
Two segment variables are assumed to be distinct if they (corresponding segments) share a crossing. For instance,  $z_a \neq z_c$, $z_b \neq z_c$, $z_d \neq z_c$, and $z_e \neq z_c$ for Figure \ref{fig:local_diagram_z}. This assumption assures that every ideal tetrahedron of $\Tcal_4$ is non-degenerate and that terms in the equation \eqref{eqn:hyp_z} are well-defined. We refer to \cite[\S 4.2]{kim2018octahedral} for details. 

\begin{figure}[!h] 
	\centering
	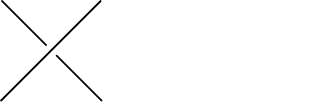
	\caption{Region variables and corners around a crossing.}
	\label{fig:crossing_w}
\end{figure}

The gluing equations of $\Tcal_5$ are expressed similarly in terms of  \emph{$w$-variables}, also called \emph{region variables}.
We have one gluing equation from each region $R$ of $D$ which is of the form
\begin{equation} \label{eqn:hyp_w}
	\tau(\kappa_1)\cdots \tau(\kappa_m) = 1.
\end{equation}
Here $\kappa_1,\ldots,\kappa_m$ are the corners of $R$ (so the region $R$ is an $m$-gon) and the $\tau$-value is given by
\begin{equation}\label{eqn:tau}
	\tau(\kappa_i) = \dfrac{ w_a w_c - w_b w_d}{( w_a -  w_d ) ( w_c -  w_d)}, \quad \tau(\kappa_j)=\dfrac{( w_b-  w_c) ( w_d- w_c)}{w_b w_d-w_a w_c} \quad 
	\textrm{for Figure~\ref{fig:crossing_w}.}
\end{equation}
Two region variables are assumed to be distinct if they (corresponding regions) share a segment. 
For instance,  $w_a,\ldots,w_d$ are pairwise distinct for Figure \ref{fig:crossing_w}.
We also assume that $w_a w_c \neq w_b w_d$ if  $w_a,\ldots,w_d$ gather around a crossing (see~Figure~\ref{fig:crossing_w} again). 
These assumptions assure that every ideal tetrahedron of $\Tcal_5$ is non-degenerate and that $\tau$-values in the equation~\eqref{eqn:tau} are well-defined. We refer to \cite[\S 4.3]{kim2018octahedral} for details.

\subsection{Ptolemy coordinates of $\Tcal_4$} \label{sec:four}

Supose that segment variables satisfying the gluing equations of $\Tcal_4$ are given with an associated developing map $\mathfrak{D} : \hat{N}\rightarrow \Hb$  and a holonomy representation $\rho:\pi_1(M)\rightarrow \psl$. Here $M$ is the knot exterior with two balls removed and $N$ is the universal cover of $M$.
For notational simplicity we identify an object in $\hat{N}$ with its image under the developing map $\mathfrak{D}$. 

\subsubsection{Ptolemy coordinates on an octahedron} \label{sec:four_octa}

\begin{figure}[!h] 
	\centering
	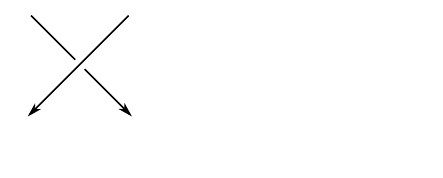
	\caption{Segment variables around a crossing.}
	\label{fig:crossing_z}
\end{figure}

\begin{figure}[!h]
	\centering
	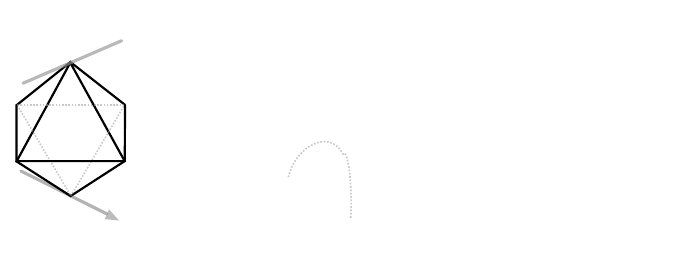
	\caption{Two developing images of $o_i$.}
	\label{fig:octa_in_H}
\end{figure}

Let us consider an ideal octahedron $o_i$ placed at a positive crossing. 
Recall that segment variables around the crossing are the coordinates of ideal vertices of $o_i$ where the bottom and top vertices are fixed by $0$ and $\infty\in\partial\Hb$, respectively. See Figures~\ref{fig:crossing_z}~(a) and~\ref{fig:octa_in_H}.

When we construct a developing map along the diagram, the octahedron $o_i$ appears twice. We denote by $O_i$ (resp., ${O}^{i}$) the one that appears when we over-pass (resp., under-pass) the crossing.
Normalizing $O_i$ so that $v_c-v_a=1$, the coordinates $v_a,v_b,v_c,v_d,v_0,$ and $v_\infty$ of the ideal vertices of $O_i$  (see Figure \ref{fig:octa_in_H}) are given by
\begin{equation}\label{eqn:coord}
	v_a=v_0+\frac{z_a}{z_c-z_a}, \ v_b=v_0+\frac{z_b}{z_c-z_a}, \ v_c=v_0+\frac{z_c}{z_c-z_a}, \ v_d=v_0+\frac{z_d}{z_c-z_a}, \ v_\infty=\infty.
\end{equation}	
Let $\Dcal : v(N)\rightarrow \psl/P$ be a  decoration described as in Proposition~\ref{prop:dev_dec}. That is,
\begin{equation}\label{eqn.rq}
	(\Dcal(v_j))(\infty)=v_j \quad \textrm{for } j \in \{a,b,c,d,0,\infty\}.
\end{equation}
There must be three $\rho$-equivariance relations on  $\Dcal(v_j)$'s,
as each ideal vertex of $\Tcal_{4}$ appears in $O_i$ twice.

It is proved in \cite[Section 5]{kim2018octahedral} that there are two loops $\gamma_1$ and $\gamma_2 \in \pi_1(M)$ whose holonomy actions are parabolic with 
\[
\rho(\gamma_1) \cdot v_a = v_c, \quad \rho(\gamma_1) \cdot v_\infty = v_\infty,\quad
\rho(\gamma_2) \cdot v_d = v_b, \quad \rho(\gamma_2) \cdot v_0 = v_0.
\]
Combining the above with the equation~(\ref{eqn:coord}), we have
\[\rho(\gamma_1)=\begin{pmatrix} 1 & 1 \\ 0 & 1 \end{pmatrix}, \quad \rho(\gamma_2) = \begin{pmatrix} 1 & v_0 \\ 0 & 1 \end{pmatrix} \begin{pmatrix} 1 & 0 \\ \Lambda & 1 \end{pmatrix} \arraycolsep=3pt\begin{pmatrix} 1 & -v_0 \\ 0 & 1 \end{pmatrix}
\]
where $\Lambda := (z_c-z_a)(\frac{1}{z_b}-\frac{1}{z_d})$.
It follows that
\begin{equation}\label{eqn:eqv1}
	\Dcal(v_c)=  \begin{pmatrix} 1 & 1 \\ 0 & 1 \end{pmatrix}\Dcal(v_a), \quad 
	\Dcal(v_b)= \begin{pmatrix} 1 & v_0 \\ 0 & 1 \end{pmatrix} \begin{pmatrix} 1 & 0 \\ \Lambda & 1 \end{pmatrix} \begin{pmatrix} 1 & -v_0 \\ 0 & 1 \end{pmatrix}\Dcal(v_d).
\end{equation}
On the other hand, as both $O_i$ and ${O}^i$ are liftings of $o_i$, there exists $\gamma_3 \in \pi_1(M)$ such that $O_i=\gamma_3 \cdot O^i$. 
The vertex ${v}^0$ of ${O}^i$ corresponding to $v_0$ coincides with the vertex $v_\infty$ of $O_i$ in $\hat{N}$ ($v_\infty$ and $v^0$ are placed at the same point $\infty$ in  Figure \ref{fig:octa_in_H}). In particular, $\Dcal({v}^0)=\Dcal(v_\infty)$ and 
\[\Dcal(v_0) = \rho(\gamma_3) \mkern 2mu \Dcal(v^0)= \rho(\gamma_3)\mkern 2mu\Dcal(v_\infty).\]
From the fact that $O_i=\gamma_3 \cdot{O}^i$ we have
\begin{equation}\label{eqn:eqv3}
\Dcal(v_0) =\begin{pmatrix} 1 & v_0 \\ 0& 1 \end{pmatrix} \arraycolsep=2pt\begin{pmatrix}0 & \sqrt{-\Lambda}^{-1} \\-\sqrt{-\Lambda} & 0 \end{pmatrix} \arraycolsep=4pt\begin{pmatrix} 1 & -v^\infty \\ 0 & 1 \end{pmatrix}\Dcal(v_\infty)
\end{equation}
where ${v}^\infty$ is the vertex of ${O}^i$ corresponding to $v_\infty$. 
Clearly, there are many choices of $\Dcal(v_j)$ that satisfies the condition~\eqref{eqn.rq}: for any $\alpha, \beta \in \Cbb$
\[\Dcal(v_j)=\arraycolsep=1.5pt\begin{pmatrix} \alpha & * \\ \beta & * \end{pmatrix}P\]
satisfies  the condition~\eqref{eqn.rq} if and only if $\alpha/\beta=v_j$. 
Among these choices, we choose 
\begin{equation} \label{eqn:deco}
\left\{
\begin{array}{l}
	\Dcal(v_a) = \dfrac{p_i}{\sqrt{z_c-z_a}} \begin{pmatrix} 1 & v_0 \\ 0 & 1 \end{pmatrix} \begin{pmatrix} z_a & * \\ z_c-z_a & * \end{pmatrix}P\\[12pt]
	\Dcal(v_d) = p_i \sqrt{z_c-z_a} \begin{pmatrix} 1 & v_0 \\ 0 & 1 \end{pmatrix} \begin{pmatrix} \frac{1}{z_c-z_a}& * \\[5pt] \frac{1}{z_d} & * \end{pmatrix}P\\[12pt]
	\Dcal(v_\infty) = \begin{pmatrix} 1 & v_0 \\ 0 & 1 \end{pmatrix} \begin{pmatrix} 1 & * \\ 0 & * \end{pmatrix}P			
\end{array} \right.
\end{equation}
for some $p_i \in \Cbb^\ast$.
Note that the above choices indeed satisfy the condition~\eqref{eqn.rq} and that we put the term $\sqrt{z_c-z_a}$ intentionally to make the resulting Ptolemy coordinates symmetric (see Figure \ref{fig:z_ptolemy_pos} below). 
It follows from the equivariance relations (\ref{eqn:eqv1})  and (\ref{eqn:eqv3}) that
\begin{equation*}
\left\{
\begin{array}{l}
	\Dcal(v_c) = \dfrac{p_i}{\sqrt{z_c-z_a}} \begin{pmatrix} 1 & v_0 \\ 0 & 1 \end{pmatrix}\begin{pmatrix}  z_c & * \\ z_c-z_a & * \end{pmatrix} P\\[12pt]
	\Dcal(v_b) = p_i \sqrt{z_c-z_a}\begin{pmatrix}1 & v_0 \\ 0 & 1 \end{pmatrix}\begin{pmatrix} \frac{1}{z_c-z_a}& * \\[5pt] \frac{1}{z_b} & * \end{pmatrix} P\\[12pt]
	\Dcal(v_0) = \begin{pmatrix} 1 & v_0 \\ 0 & 1 \end{pmatrix}\begin{pmatrix}0 & * \\ -\sqrt{-\Lambda} & * \end{pmatrix}P	
\end{array} \right.
\end{equation*}

\begin{figure}[!h] 
\centering  \small
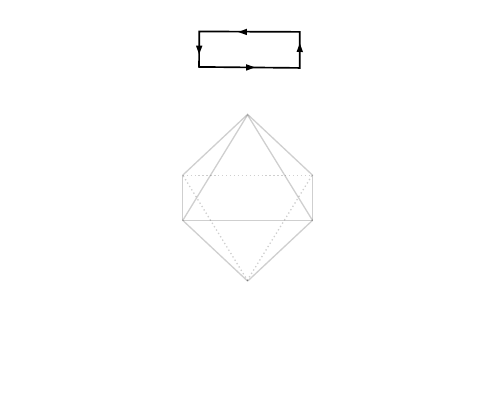
\caption{Short-edges of a truncated octahedron.}
\label{fig:z_short_edge}
\end{figure}

We finally compute Ptolemy coordinates on $o_i$ by using the formula~\eqref{eqn:detf}.  
The result is given in Figure~\ref{fig:z_ptolemy_pos}.
\begin{figure}[!h]
\centering
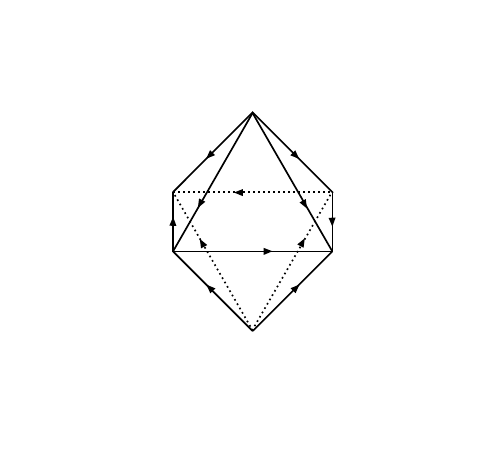
\caption{Ptolemy coordinates for a positive crossing.}
\label{fig:z_ptolemy_pos}
\end{figure}
Note that we have not used the fact that the holonomy representation $\rho$ is a $\psl$-representation; all  computations so far work for $\sl$. Indeed, Ptolemy coordinates computed in Figure~\ref{fig:z_ptolemy_pos} satisfy the equation~\eqref{eqn:ptolemy_obs} for the trivial obstruction cocycle $\sigma$. It follows that we obtain a natural $(\sl,P$)-cocycle
\[\Bcal_i : e(\overline{o}_i) \rightarrow \sl\]
where $\overline{o}_i$ is a truncated octaehdron of $o_i$.
We present some short-edge parameters of $\Bcal_i$ here:
\begin{equation} \label{eqn:z_short} 
\begin{array}{llll}
	\cl{a}= \frac{z_a-z_d}{z_c-z_a} & \cl{b} = \frac{z_b-z_a}{z_c-z_a}&\cl{c} = 	\frac{z_c-z_b}{z_c-z_a}&  \cl{d} =	\frac{z_d-z_c}{z_c-z_a}\\[6pt]
	\cl{e} = \frac{z_b}{p_i^2 z_a(z_b-z_a)} & \cl{f} = \frac{z_d}{p_i^2 z_a(z_a-z_d)}& \cl{g} = 	\frac{1}{p_i^2(z_d-z_a)}&  \cl{h} =	\frac{1}{p_i^2(z_a-z_b)}\\[6pt]
	\cl{i} = 	\frac{z_a z_b}{p_i^2(z_a-z_b)}& \cl{j} = 	\frac{z_bz_c}{p_i^2(z_b-z_c)}& \cl{k} = 	\frac{z_b^2}{p_i^2 (z_c-z_b)}&  \cl{l} = 	\frac{z_b^2}{p_i^2(z_b-z_a)} \\[6pt]
	\cl{m} = 	\frac{z_d(z_c-z_b)}{z_c(z_b-z_d)}& \cl{n} = 	\frac{z_b(z_d-z_c)}{z_c(z_b-z_d)}& \cl{o} =	\frac{z_b(z_a-z_d)}{z_a(z_b-z_d)}&  \cl{p} = 	\frac{z_d(z_b-z_a)}{z_a(z_b-z_d)} \\[6pt]
	\cl{q} = 	\frac{z_b}{p_i^2z_c(z_b-z_c)}& \cl{r} = 	\frac{z_d}{p_i^2z_c(z_c-z_d)}& \cl{s} = 	\frac{1}{p_i^2(z_d-z_c)}& \cl{t} = 	\frac{1}{p_i^2 (z_c-z_b)} \\[6pt]
	\cl{u} = 	\frac{z_a z_d}{p_i^2 (z_a-z_d)}&  \cl{v} = 	\frac{z_c z_d}{p_i^2(z_d-z_c)}& \cl{w} = 	\frac{z_d^2}{p_i^2(z_c-z_d)}& \cl{x} = 	\frac{z_d^2}{p_i^2(z_d-z_a)} 
\end{array}
\end{equation}
where $\cl{a},\ldots,\cl{x}$ are short-edges of $\overline{o}_i$ as in Figure ~\ref{fig:z_short_edge}.

We compute Ptolemy coordinates for a negative crossing similarly. It turns out that a negative crossing as in Figure~\ref{fig:crossing_z}~(b) has the same Ptolemy coordinates as in Figure \ref{fig:z_ptolemy_pos} except that the term $\sqrt{z_c-z_a}$ is replaced by $\sqrt{z_a-z_c}$.  
Also, short-edge parameters of $\Bcal_i$ for a negative crossing are the same as those in the equation~\eqref{eqn:z_short}  except that the sign for $\cl{a},\ldots,\cl{d}$ changes.

\subsubsection{Gluing octahedra} \label{sec:four_comp}
As the ideal octahedra $o_1,\ldots,o_n$ are glued along their faces (to form the knot complement minus two points), we require certain relations on the variables $p_1,\ldots,p_n$ to glue the natural cocycles $\Bcal_i : e(\overline{o}_i) \rightarrow \sl$ compatibly. See \cite[\S 3.1]{kim2018octahedral} for how the octahedra $o_1,\ldots,o_n$ are glued.
One easily checks that the relations on $p_1,\ldots,p_n$ consist of 
\begin{equation}\label{eqn:p_relation}
	\begin{cases}
	p_i \sqrt{1- \dfrac{z_b}{z_c}}= p_j \sqrt{1- \dfrac{z_c}{z_e}} & \textnormal{for Figure \ref{fig:local_diagram_z}~(a)} \\[10pt]
	p_i \sqrt{1- \dfrac{z_c}{z_b}}=p_j \sqrt{1- \dfrac{z_e}{z_c}} & \textnormal{for Figure \ref{fig:local_diagram_z}~(b)} \\[10pt]
	p_i \sqrt{1- \dfrac{z_b}{z_c}}=p_j  \sqrt{1- \dfrac{z_e}{z_c}} & \textnormal{for Figure \ref{fig:local_diagram_z}~(c)}  \\[10pt]
	p_i  \sqrt{1- \dfrac{z_c}{z_b}}=p_j  \sqrt{1- \dfrac{z_c}{z_e}} & \textnormal{for Figure \ref{fig:local_diagram_z}~(d)} 
\end{cases}
\end{equation}
for each segment of $D$ where $i$ and $j$ are the indices of ideal octahedra placed at the crossing in the left and right of the segment, respectively. 
It turns out that it may not be possible to glue all $\Bcal_i$ compatibly. However, it is possible if we project them down to $\psl$ (see Proposition~\ref{prop:compatible} below).
Since $\Bcal_i$ represents the same natural $(\psl,P)$-cocycle for both $p_i$ and $-p_i$, we may regard $p_i$ as an element of $\Cbb^\ast/\pmo$.

\begin{proposition} \label{prop:compatible}
	Regarding each $p_i$ as an element of $\Cbb^\ast / \pmo$, there exists a collection of $p_1,\ldots,p_n $ satisfying the equation \eqref{eqn:p_relation} (up to sign) for all segments of $D$. Moreover, such a collection is unique up to scaling.
\end{proposition}
\begin{proof}
	It is clear that all $p_i$'s are determined  successively along the diagram (due to the equation \eqref{eqn:p_relation}) if one of them is given. Thus the uniqueness part is obvious. 
	To prove the existence part, we consider a region $R$ of $D$ as in Figure \ref{fig:region}. 
	Rewriting the equation \eqref{eqn:p_relation} for Figure \ref{fig:local_diagram_z}~(a) by using the equation \eqref{eqn:hyp_z}, we have (up to sign)
	\[p_i \sqrt{1-\dfrac{z_a}{z_c}} = p_j \sqrt{1-\dfrac{z_c}{z_d}}\, .\]
	It follows that	relations on $p_1,\ldots,p_m$ arising from the boundary segments of $R$ are
	\[p_i \sqrt{1-\dfrac{z_i}{z_{i+1}}} = p_{i+1} \sqrt{1-\dfrac{z_{i+1}}{z_{i+2}}}, \quad 1 \leq i\leq m.\]
	Here the index is taken  modulo $m$. These relations are compatible, in the sense that their product along the boundary of $R$ results in the trivial identity. 
	We have proved the case that the boundary of $R$ is alternating as in Figure \ref{fig:region}, and one can prove other cases similarly. This proves that all the relations on $p_1,\ldots,p_n$ are compatible.
	\begin{figure}[!h]
		\centering
		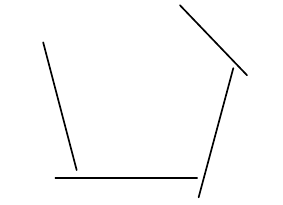
		\caption{Segment variables around a region.}
		\label{fig:region}
	\end{figure}
\end{proof}

\begin{remark}
	As we will see in Section \ref{sec:obs}, exact values of
	$p_1,\ldots,p_n$ in Proposition~\ref{prop:compatible} are not that important. They vanish eventually in the computation for a holonomy representation (see Section~\ref{sec:simplify}).
\end{remark} 
Restating Proposition \ref{prop:compatible}, we obtain:
\begin{theorem} \label{thm:compatible}
	The natural cocycles $\Bcal_i : e(\overline{o}_i) \rightarrow \sl$ $(1 \leq i\leq n)$ are well-glued if we project them down to $\psl$.
	In particular, they form a natural cocyle $\Bcal:e(M)\rightarrow \psl$.
\end{theorem}

As we explained in Section~\ref{sec:ptol}, the natural cocyle $\Bcal:e(M)\rightarrow \psl$ defines a map 
\[c:e(\Tcal_4)\rightarrow \Cbb^\ast\]
after some sign-choices.
Also, Proposition~\ref{prop:ptolob} shows that $c\in P^\sigma(\Tcal_4)$ for some $\sigma \in Z^2(\Tcal_4;\pmo)$ whose class in $H^2(M, \partial M ;\pmo)$ is the obstruction class of the holonomy representation $\rho$. We analyze this class in Section \ref{sec:obs}.

\subsubsection{Simplification} \label{sec:simplify}
\begin{figure}[!h] 
	\centering  	
	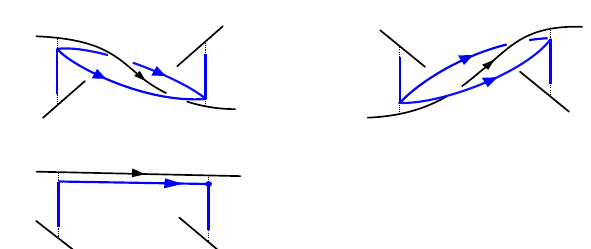
	\caption{Local configurations of $G$ for Figure \ref{fig:local_diagram_z}.}
	\label{fig:cellcomp}
\end{figure}

The Ptolemy coordinates of $\Tcal_{4}$ that we have computed  seems somewhat complicated, as there are too many edges. 
Motivated by \cite{thistlethwaite2014alternative} and \cite{neumann2016intercusp}, we simplify them by considering a certain graph $G$ defined as follows.
\begin{enumerate}
	\item  Place a \emph{vertical edge} at each crossing of $D$ so that it  represents the central axis of the ideal octahedron $o_i$.
	\item For each segment of $D$ join two vertical edges by either one or two \emph{horizontal edges} along the segment as in Figure \ref{fig:cellcomp}. 
\end{enumerate}
Note that any loop in $M$ can be homotoped to an edge-path of $G$.

Each vertical edge of $G$ is a long-edge of $M$ and each horizontal edge of $G$ is the composition of short-edges of $M$. Thus, the natural cocycle $\Bcal:e(M)\rightarrow \psl$ in Theorem~\ref{thm:compatible} endows each edge of $G$ with a \emph{long} or \emph{short-edge parameter} according to the edge type. Precisely, the long-edge parameter $c(e) \in \Cbb^\ast / \pmo$ of a vertical edge $e$ is given by
\begin{equation*}
	c(e)=
	\begin{cases}
			\begin{array}{cc}
					\sqrt{(z_a-z_c)\left(\frac{1}{z_d}-\frac{1}{z_b}\right)} & \textrm{for  Figure \ref{fig:crossing_z}(a)} \\[12pt]
					\sqrt{(z_a-z_c)\left(\frac{1}{z_b}-\frac{1}{z_d}\right)} & \textrm{for  Figure \ref{fig:crossing_z}(b)}
				\end{array}		
	\end{cases}
\end{equation*} and the short-edge parameter $c(e)\in\Cbb$ of a  horizontal edge $e$ is given by
\begin{equation*}
	c(e)=\left\{
	\begin{array}{cc}
		\dfrac{z_c}{z_a-z_b} - \dfrac{z_d}{z_d-z_e} & \textrm{for  Figure \ref{fig:cellcomp}(a)} \\[12pt]
		\dfrac{z_a}{z_a-z_b} - \dfrac{z_c}{z_d-z_e} & \textrm{for  Figure \ref{fig:cellcomp}(b)} \\[12pt]
		\dfrac{z_c}{z_a-z_b} - \dfrac{z_c}{z_d-z_e} & \textrm{for  Figure \ref{fig:cellcomp}(c)} \\[12pt]
		\dfrac{z_a}{z_a-z_c} - \dfrac{z_d}{z_d-z_e} & \textrm{for  Figure \ref{fig:cellcomp}(d)}
	\end{array}
	\right.
\end{equation*} Note that the variable $p_i$ does not appear in the above equations.

\begin{remark} In  \cite{neumann2016intercusp} (resp., \cite{thistlethwaite2014alternative}) long-edge and short-edge parameters of $G$ are called intercusp and translation parameters (resp., crossing and edge labels).
\end{remark}

\subsection{Ptolemy coordinates of $\Tcal_{5}$}\label{sec:five}
Suppose that region variables  satisfying the gluing equations of $\Tcal_5$ are given with a  holonomy representation $\rho:\pi_1(M)\rightarrow \psl$.
In this section, we compute Ptolemy coordinates of $\Tcal_5$ in terms of region variables by using the same strategy that we used in Section \ref{sec:four}. Most of the computations in this section will be similar to those in Section \ref{sec:four}, hence some of them will be omitted to avoid tedious repetitions.

\subsubsection{Ptolemy coordinates on an octahedron}

\begin{figure}[!h] 
	\centering
	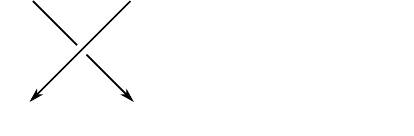
	\caption{Region variables around a crossing.}
	\label{fig:crossing_ww}
\end{figure}

Let us consider an ideal octahedron $o_i$ placed at a positive crossing. We denote region variables around the crossing by $w_a, \ldots,w_d$  as in Figure \ref{fig:crossing_ww}~(a). 
Ratios of these variables are cross-ratios of ideal tetrahedra of $o_i$ (see \cite[\S 4.3]{kim2018octahedral}):
\begin{equation} \label{eqn:cross_ratio}
	\left[v_b,v_\infty,v_c,v_a\right]=\frac{w_a}{w_b},\ \left[v_d,v_\infty,v_a,v_c\right]=\frac{w_c}{w_d},\
	\left[v_0,v_c,v_d,v_b\right] =\frac{w_c}{w_b}
\end{equation} 
where $v_a,v_b,v_c,v_d,v_0,v_\infty$  are the vertices of  the ideal octahedron $O_i$ as in Figure~\ref{fig:octa_in_H}.
Here we use the convention
\[[x,y,z,w]= \frac{(x-w)(y-z)}{(x-z)(y-w)}. \]
Recall that $O_i$ is normalized so that $v_c-v_a=1$ and $v_\infty=\infty$.  We may assume further that $v_a=0$ and $v_c=1$, as we already saw in Section \ref{sec:four} that only relative coordinates of the ideal vertices matter when we compute Ptolemy coordinates. A straightforward computation using the equation \eqref{eqn:cross_ratio} gives 
\begin{equation*}
	v_b= \dfrac{w_a}{w_a-w_b},\ v_d=\dfrac{w_d}{w_d-w_c},\ v_0=\dfrac{w_a-w_d}{w_a-w_b+w_c-w_d}.
\end{equation*} 
Note that we allow the case $v_0=\infty$, that is, $w_a-w_b+w_c-w_d=0$. 

\begin{figure}[!h]
	\centering
	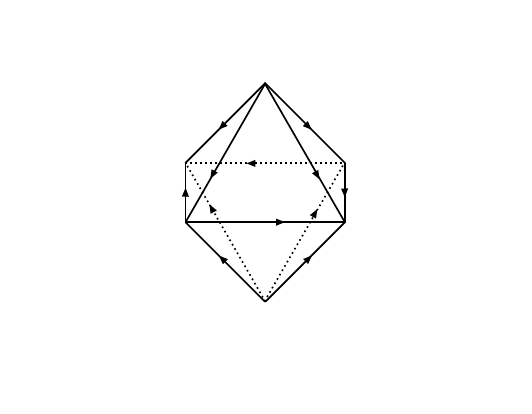
	\caption{Ptolemy coordinates for a positive crossing.}
	\label{fig:w_ptolemy_pos}
\end{figure}

We use the same loops $\gamma_1,\gamma_2,\gamma_3 \in \pi_1(M)$ in Section \ref{sec:four} to compute the $\rho$-equivariance relations on decoration $\Dcal$.
If $v_0 \neq \infty$, the computation would be exactly the same as in Section~\ref{sec:four}. The $\rho$-equivariance relations result in the same relations~\eqref{eqn:eqv1} and \eqref{eqn:eqv3} except that $\Lambda$ is now expressed as
\begin{equation} \label{eqn:Lambda}
	\Lambda := \frac{(w_a-w_b+w_c-w_d)^2}{w_a w_c - w_b w_d}.
\end{equation}
Letting
\begin{equation} \label{eqn:deco2}
	\Dcal(v_a) =q_i \begin{pmatrix} 0 & * \\ 1 & * \end{pmatrix}P,\quad
	\Dcal(v_d) = \dfrac{1}{q_i} \begin{pmatrix} w_d& * \\ w_d-w_c & * \end{pmatrix}P,\quad 
	\Dcal(v_\infty) = \begin{pmatrix} 1 & * \\ 0 & * \end{pmatrix} P
\end{equation}
for some $q_i\in\Cbb^\ast$, we have
\begin{equation*} 
	\Dcal(v_c) = q_i \begin{pmatrix} 1 & * \\ 1 & * \end{pmatrix}P, \quad 
	\Dcal(v_b) = \dfrac{1}{q_i} \begin{pmatrix}w_a & * \\ w_a-w_b & * \end{pmatrix}P ,\quad
	\Dcal(v_0) = \sqrt{-\Lambda}\begin{pmatrix}v_0& * \\ 1 & * \end{pmatrix}P
\end{equation*} 
from the equivariance relations.
Then one can compute Ptolemy coordinates on the ideal octahedron $o_i$ by using the equation \eqref{eqn:detf}. 
See Figure \ref{fig:w_ptolemy_pos} for the result.
If $v_0=\infty$, then a holonomy representation for $\gamma_2$ and $\gamma_3$ should change to
\[\rho(\gamma_2) = \begin{pmatrix} 1 & v_b-v_d \\ 0 & 1 \end{pmatrix},\quad \rho(\gamma_3)=\begin{pmatrix} \frac{\sqrt{w_b w_d-w_aw_c}}{w_a-w_b} & \frac{-w_d}{\sqrt{w_b w_d-w_aw_c}} \\[5pt] 0 &  \frac{w_a-w_b}{\sqrt{w_b w_d-w_aw_c}} \end{pmatrix}.\]
However, it turns out that Ptolemy coordinates on $o_i$ remain the same as in Figure~\ref{fig:w_ptolemy_pos}.

As the above computation works in $\sl$,
the Ptolemy coordinates in Figure~\ref{fig:w_ptolemy_pos} satisfy the equation~\eqref{eqn:ptolemy_obs} for the trivial obstruction cocycle $\sigma$, and induce a natural cocycle 
\[\Bcal_i : e(\overline{o}_i) \rightarrow \sl\]
on the truncated octahedron $\overline{o}_i$. 
We present some short-edge parameters of $\Bcal_i$ here.
\begin{equation} \label{eqn:w_short} 
	\begin{array}{llll}
		\cl{a} = \frac{w_d}{w_c-w_d} & \cl{b} = \frac{w_a}{w_a-w_b}& 
		\cl{c} = \frac{w_b}{w_b-w_a}&  \cl{d} = \frac{w_c}{w_d-w_c}\\[6pt]
		\cl{e} = \frac{w_bw_d-w_aw_c}{q_i^2 w_a(w_a-w_d)}& \cl{f} = 	\frac{w_bw_d-w_aw_c}{q_i^2w_d(w_d-w_a)}&  \cl{g} = \frac{w_d-w_c}{q_i^2w_d}& \cl{h} = 	\frac{w_b-w_a}{q_i^2 w_a}\\[6pt]
		\cl{i} = \frac{q_i^2(w_d-w_a)}{w_a(w_bw_d-w_aw_c)}& \cl{j} = 	\frac{q_i^2(w_b-w_c)}{w_b(w_bw_d-w_aw_c)}& \cl{k} = 	\frac{q_i^2}{w_b(w_b-w_a)}&  \cl{l} = \frac{q_i^2}{w_a(w_a-w_b)} \\[6pt]
		\cl{m} = 	\frac{w_b}{w_c-w_b}& \cl{n} = 	\frac{w_c}{w_b-w_c}& \cl{o} =	\frac{w_d}{w_d-w_a}&  \cl{p} = \frac{w_a}{w_a-w_d} \\[6pt]
		\cl{q} = 	\frac{w_bw_d-w_aw_c}{q_i^2w_b(w_b-w_c)}&\cl{r} = 	\frac{w_bw_d-w_aw_c}{q_i^2w_c(w_c-w_b)}& \cl{s} = 	\frac{w_d-w_c}{q_i^2 w_c}& \cl{t} = 	\frac{w_b-w_a}{q_i^2 w_b} \\[6pt]
		\cl{u} = 	\frac{q_i^2(w_d-w_a)}{w_d(w_bw_d-w_aw_c)}&  \cl{v} = 	\frac{q_i^2(w_b-w_c)}{w_c(w_bw_d-w_aw_c)}&	\cl{w} = 	\frac{q_i^2}{w_c(w_c-w_d)}& \cl{x} = 	\frac{q_i^2}{w_d(w_d-w_c)}
	\end{array}
\end{equation}
where $\cl{a},\ldots,\cl{x}$ are short-edges of $\overline{o}_i$ as in Figure ~\ref{fig:z_short_edge}.

\begin{figure}[!h]
	\centering
	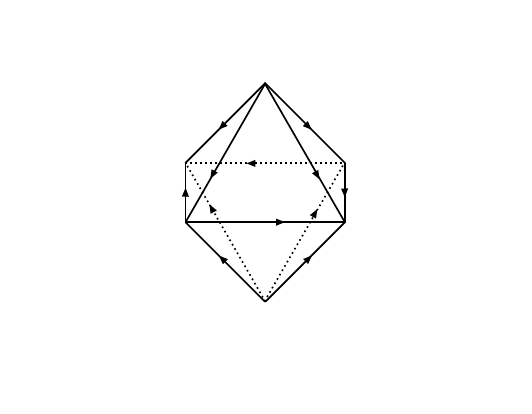
	\caption{Ptolemy coordinates for a negative crossing.}
	\label{fig:w_ptolemy_neg}
\end{figure}

We similarly compute Ptolemy coordinates for a negative crossing as in Figure~\ref{fig:crossing_ww}~(b). See Figure~\ref{fig:w_ptolemy_neg} for the result.
Short-edge parameters of $\Bcal_i$ for a negative crossing are the same  as  in the equation~\eqref{eqn:w_short} except that the sign for $\cl{a},\ldots,\cl{l}$ and $\cl{q} ,\ldots,\cl{x}$ changes.

\subsubsection{Gluing octahedra}

\begin{figure}[!h] 
	\centering  	
	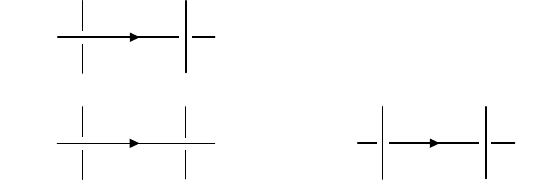
	\caption{Region variables around a segment.}
	\label{fig:local_diagram_w}
\end{figure}

As the ideal octahedra $o_1,\ldots,o_n$ are glued along their faces, we require to certain relations on the variables $q_1,\ldots,q_n$ to glue the natural cocycles $\Bcal_i : e(\overline{o}_i) \rightarrow \sl$ compatibly.  
For convenience reasons, we define an \emph{$\eta$-value at a crossing} by $\eta:=w_bw_d-w_aw_c$ for Figure \ref{fig:crossing_ww}~(a) and  $\eta:=w_aw_c-w_bw_d$ for Figure \ref{fig:crossing_ww}~(b).
Then one can check that the relations on  $q_1,\ldots,q_n$ consist of 
\begin{equation}\label{eqn:w_p_relation}
	\dfrac{q_i}{q_j}=
	\left \{
	\begin{array}{ll}
		\dfrac{w_c-w_d}{\sqrt{\eta_j}}& \textnormal{for Figure \ref{fig:local_diagram_w}~(a)} \\[12pt]
		\dfrac{\sqrt{\eta_i}}{w_c-w_d}    & \textnormal{for Figure \ref{fig:local_diagram_w}~(b)}  \\[12pt]
		1 & \textnormal{for Figure \ref{fig:local_diagram_w}~(c)} \\[5pt]
		\dfrac{\sqrt{\eta_i}}{\sqrt{\eta_j}} & \textnormal{for Figure \ref{fig:local_diagram_w}~(d)}
	\end{array}
	\right.
\end{equation}
for each segment of $D$ where $i$ and $j$ are the indices of ideal octahedra placed at the crossing in the left and right of the segment, respectively. 
It turns out that it may not be possible to  glue all $\Bcal_i$'s compatibly. However, it is possible if we project them down to $\psl$ (see Proposition \ref{prop:compatible_w} below). Since  $\Bcal_i$ represents the same natural $(\psl,P)$-cocycle for both $q_i$ and $-q_i$, we may regard $q_i$ as an element of $\Cbb^\ast/\pmo$.

\begin{proposition} \label{prop:compatible_w}
	Regarding each $q_i$ as an element of $\Cbb^\ast / \pmo$, there exists a collection of $q_1,\ldots,q_n $ satisfying the equation \eqref{eqn:w_p_relation} (up to sign) for all segments of $D$. Moreover, such a collection is unique up to scaling.
\end{proposition}
\begin{proof}
	As in the proof of Proposition \ref{prop:compatible}, we consider a region  $R$ with the equations \eqref{eqn:w_p_relation} arising from its boundary segments. It suffices to check that the product of these equations results in the trivial identity up to sign.
	
	Suppose that $R$ is given as in Figure \ref{fig:region_w}. Then relations on $q_1,\ldots,q_m$ arising from its boundary segments are 
	\begin{equation} \label{eqn:comp_w}
		\dfrac{q_i}{q_{i+1}} =  \dfrac{w_{2i+1} - w_0}{\sqrt{\eta_{i+1}}}, \quad 1 \leq i \leq m
	\end{equation}
	where $\eta_i = w_{2i-i}w_{2i} - w_0 w_{2i}$. Here the index of $q_i$ (resp., $w_i$) is taken in modulo $m$ (resp., $2m$). 
	On the other hand, the equation \eqref{eqn:hyp_w} for the region $R$ is
	\[1=\prod_{i=1}^m \dfrac{(w_{2i-1}-w_0)(w_{2i+1}-w_0)}{w_{2i-1}w_{2i+1}-w_0 w_{2m}}= \prod_{i=1}^m \dfrac{(w_{2i+1}-w_0)^2}{\eta_i}.\]
	Therefore, the product of the equations \eqref{eqn:comp_w} over $1 \leq i\leq n$ is the trivial identity (up to sign). If we change the orientation of a boundary segment of $R$, the sign of two etas should change. Therefore, the product of the equations \eqref{eqn:comp_w} over $1 \leq i\leq n$ does not depend on the orientation of segments.
	We have proved the case when the boundary of $R$ is alternating, and one can prove other cases similarly.
	\begin{figure}[!h]
		\centering
		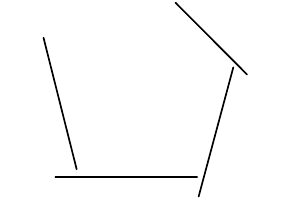
		\caption{Region variables around a region.}
		\label{fig:region_w}
	\end{figure}	
\end{proof}
Restating Proposition \ref{prop:compatible_w}, we obtain:

\begin{theorem} \label{thm:compatible_w}
	The natural cocycles $\Bcal_i : e(\overline{o_i}) \rightarrow \sl$ $(1 \leq i\leq n)$ are well-glued if we project them down to $\psl$.
	In particular, they form a natural cocyle $\Bcal:e(M)\rightarrow \psl$.
\end{theorem}
As we explained in Section~\ref{sec:ptol}, the natural cocyle $\Bcal:e(M)\rightarrow \psl$ defines a map 
\[c:e(\Tcal_5)\rightarrow \Cbb^\ast\]
after some sign-choices, where $c\in P^\sigma(\Tcal_5)$ for some $\sigma \in Z^2(\Tcal_5;\pmo)$ whose class in $H^2(M, \partial M ;\pmo)$ is the obstruction class of the holonomy representation $\rho$. 
We analyze this class in Section \ref{sec:obs}.

\subsubsection{Simplification}

\begin{figure}[!h] 
	\begin{subfigure}[!h]{\textwidth}	
		\centering
		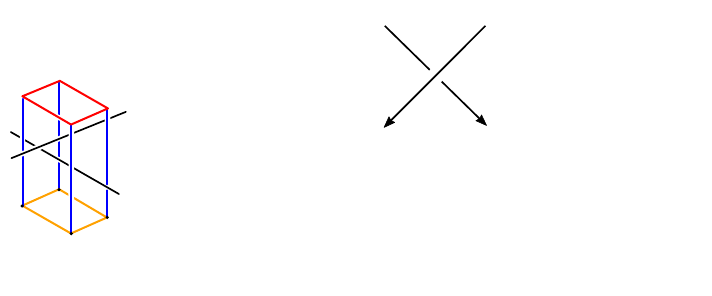
		\caption{Positive crossing.}
	\end{subfigure}
	
	\begin{subfigure}[!h]{\textwidth}	
		\centering
		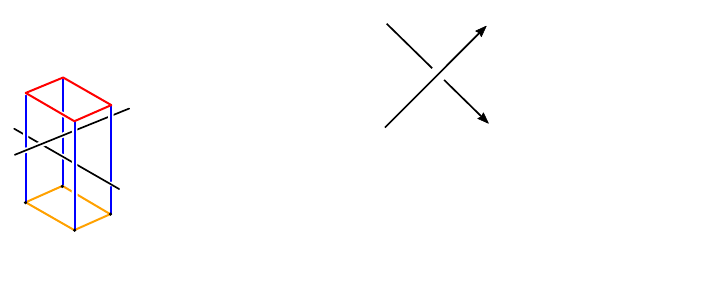
		\caption{Negative crossing.}
	\end{subfigure}
	\caption{Short-edge parameters of $X$.}
	\label{fig:cellcomp2}
	
\end{figure}

We simplify the Ptolemy coordinates of $\Tcal_5$ by considering a graph $X$ defined as follows.  
\begin{enumerate}
	\item For each region of $D$ place an edge  penetrating the region vertically ($X$ then has $N+2$ $1$-cells with $2N+4$ $0$-cells); 
	\item For each pair of adjacent regions of $D$ join the top and bottom vertices of the penetrating edges respectively (see Figure~\ref{fig:cellcomp2}).
\end{enumerate} 
Note that any loop in $M$ can be homotoped to an edge-path of $G$.

Each vertical edge of $X$ is a long-edge of $M$ and each horizontal edge of $X$ is the composition of short-edges of $M$. Thus, the natural cocycle $\Bcal:e(M)\rightarrow \psl$ in Theorem~\ref{thm:compatible_w} endows each edge of $G$ with a \emph{long} or \emph{short-edge parameter} according to the edge type.  One easily checks from Figure~\ref{fig:w_ptolemy_pos} or ~\ref{fig:w_ptolemy_neg} that the long-edge parameter of a vertical edge of $X$ is exactly the region variable.
The short-edge parameters of horizontal edges of $X$ are given in Figure~\ref{fig:cellcomp2}.

\section{The obstruction class and the cusp shape} \label{sec:obs}

Let $M$ be the knot exterior of $K \subset S^3$ with two balls removed and $D$ be a diagram of $K$ with $n$ crossings $c_1,\ldots,c_n$.
If segment variables are given, we define two values $\sigma(c_i)$ and $\lambda(c_i)$ at each crossing $c_i$ as
\[\sigma(c_i) := \dfrac{z_a}{z_c}, \quad \lambda(c_i) := \dfrac{z_b z_d \mkern1mu(z_c-z_a)}{z_a z_c\mkern1mu(z_b-z_d)}\] where $z_a,\ldots,z_d$ are  segment variables around the crossing as in Figure \ref{fig:crossing_z}~(a) or~(b).   
Similarly, if region variables are given, we define two values $\sigma(c_i)$ and $\lambda(c_i)$ at each crossing $c_i$ as
\[\sigma(c_i)  :=\dfrac{w_a-w_d}{w_b-w_c}, \quad \lambda(c_i):= \dfrac{w_a w_c - w_b w_d}{(w_a-w_d)(w_c-w_b)}\]
where $w_a,\ldots, w_d$ are region variables around the crossing as in Figure \ref{fig:crossing_ww}~(a) or~(b).

In this setting, we can write down the obstruction class and the cusp shape of $\rho$ explicitly in terms of $z$- or $w$-variables as follows.

% Then the obstruction class viewed as an element of $\pmo$ is computed as follows.
\begin{theorem} \label{thm:obs}
	Suppose that segment variables (resp., region variables) satisfying the gluing equations of $\Tcal_4$ (resp., $\Tcal_5$) are given with a holonomy representation $\rho:\pi_1(M)\rightarrow \psl$. Then we have
	\begin{equation} \label{eqn:pm1}
		\sigma:=\prod_{i=1}^n \sigma(c_i) \in \pmo.
	\end{equation}
	Moreover, $\sigma$ agrees with the obstruction class of  $\rho$ under the duality 
	\[H^2(M,\partial M;\pmo) \simeq H_1(M;\pmo) \simeq \pmo.\] 
\end{theorem}

\begin{theorem} \label{thm:longitude} 
	Suppose that segment variables (resp., region variables) satisfying the gluing equation of $\Tcal_4$ (resp., $\Tcal_5$) are given with a holonomy representation $\rho:\pi_1(M)\rightarrow \psl$. Then the cusp shape of $\rho$ is given by 
	\[ \sum_{i=1}^{n} \lambda(c_i)-w(D)\]
	where $w(D)$ is the writhe of $D$.
\end{theorem}	

Before we prove the above theorems, we present how they work in examples with some remarks.

%We remark that the equation \eqref{eqn:pm1} itself is rather unexpected from the gluing equations of $\Tcal_4$ or $\Tcal_5$.

\begin{figure}[!h]
	\centering
	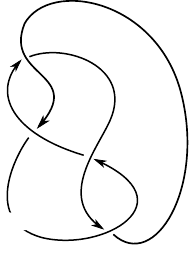
	\caption{The figure-eight knot. }
	\label{fig:figure_one}
\end{figure}

\begin{example}\label{ex:figure_obs}
	Let us consider the figure-eight knot with a diagram given as in Figure~\ref{fig:figure_one}.
	It is proved in \cite{kim2018octahedral} that segment variables 
	\[(z_1,\ldots,z_8) = \left( 1+i,\ i(1+\sqrt{3}),\ \frac{-1+i\sqrt{3}}{-1+\sqrt{3}},\ 2i,\ -1+i,\ i,\ \frac{1+i\sqrt{3}}{-1+\sqrt{3}},\ -\frac{2i}{-1+\sqrt{3}} \right)\]
	satisfy the gluing equations of $\Tcal_4$ (these are obtained by putting $p=1+i,\ q=1-i,\ r=1$ in Example 4.6 of \cite{kim2018octahedral}). Applying Theorem \ref{thm:obs}, the obstruction class of the holonomy representation is given by
	\[\sigma=\prod_{1\leq i \leq 4} \sigma(c_i)  = \frac{1+i}{i(1+\sqrt{3})}\cdot \frac{-1+i\sqrt{3}}{2i(-1+\sqrt{3})} \cdot \frac{-1+i}{i} \cdot \frac{1+i\sqrt{3}}{-2i}= -1.\]
	Applying Theorem~\ref{thm:longitude}, the cusp shape of the holonomy representation $\rho$ is given by
	\[\sum_{1 \leq i \leq 4} \lambda(c_i) - w(D) = \left(1+i(\sqrt{3}-2)+\frac{1+i\sqrt{3}}{2+2\sqrt{3}}+\frac{1+i}{\sqrt{3}+i} +\frac{2+\sqrt{3}+i}{-2+\sqrt{3}-i} \right)-0 = 2\sqrt{3}i.\]
\end{example}

\begin{figure}[!h]
	\centering
	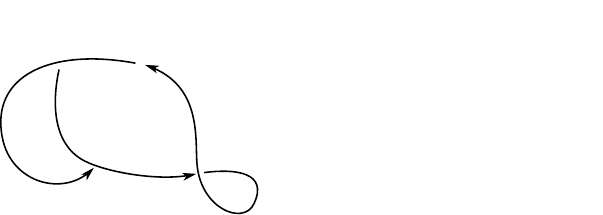
	\caption{The trefoil knot with a kink.}
	\label{fig:trefoil_one}
\end{figure}

\begin{example}
	\label{ex:trefoil_obs} 
	Let us consider the trefoil knot with a diagram given as in Figure~\ref{fig:trefoil_one}. It is proved in \cite{cho2016existence} that region variables $(w_1,\dots,w_6)=(5,3,7,2,1,8)$ satisfy the gluing equations of $\Tcal_5$. Applying Theorem \ref{thm:obs}, the obstruction class of the holonomy representation is given by
	\[\sigma = \prod_{1 \leq i \leq 4} \sigma(c_i)  = -2\cdot\frac{1}{3} \cdot 1 \cdot \dfrac{3}{2} = -1.\]
	Applying Theorem~\ref{thm:longitude}, the cusp shape of the holonomy representation $\rho$ is given by
	\[\sum_{1 \leq i \leq 4} \lambda(c_i) - w(D) = \left(-\frac{9}{2}+\frac{4}{3}-1 +\frac{1}{6} \right)-2 =-6.\]
	
\end{example}

\begin{remark} If a diagram has a kink, then  the ideal triangulation $\Tcal_4$ has a degenerate edge. In particular, the gluing equations of $\Tcal_4$ has no solution (see \cite{segerman2011pseudo}). 
	Meanwhile, it is proved in \cite{cho2016existence} that the gluing equations  of $\Tcal_5$ produce all boundary-parabolic representations except the trivial one for any knot diagram.
\end{remark}

\begin{remark}
	The geometric representation of a hyperbolic knot has a non-trivial obstruction class (see e.g. \cite{calegari2006real}). Thus a simple corollary of Theorem~\ref{thm:obs} is: if $\sigma \neq -1$, then the holonomy representation is not geometric.
\end{remark} 

\subsection{Proof of Theorem \ref{thm:obs}}\label{sec:proof}

Let  $\Bcal:e(M)\rightarrow \psl$ be the natural cocycle given by Theorem~\ref{thm:compatible} or \ref{thm:compatible_w} 
and $c':e(\Tcal)\rightarrow \Cbb^\ast/\pmo$, where $\Tcal=\Tcal_4$ or $\Tcal_5$, be 
the corresponding long-edge parameters with sign-ambiguity. 

The manifold $M$ has two  boundary spheres and one boundary torus. 
It is proved in \cite{kim2018octahedral} 
that 
the edges of $\Tcal$ that join the boundary torus with one of the boundary spheres are in one-to-one correspondence with the arcs of $D$.
We denote by $a_i$ ($1 \leq i \leq n$) an arc of $D$ and $v_i$ the corresponding edge of $\Tcal$. We may assume that $a_{i-1}$ and $a_{i}$ are consecutive along the diagram and denote by $c_i$ the crossing lying between $a_{i-1}$ and $a_i$. See Figure~\ref{fig:main_proofs}. Here the index is taken in modulo $n$.
%	Assuming that the sphere lies above the diagram
%	One of them lies over the diagram while the other one lies below. We may assume that the vertex $p$ lies above the diagram. 
%	The edges of $M$ joining $p$ with the ideal vertex corresponding to the knot $K$ are called \emph{over-edges}  and they are in one-to-one correspondence with the (over-)arcs of $D$. We enumerate the arcs by $a_i$ ($1 \leq i\leq n$) along the diagram and denote by $v_i$  the over-edge corresponding to $a_i$. Also, we may assume that a crossing $c_i$ lies between $a_{i-1}$ and $a_{i}$ (see Figure \ref{fig:main_proofs}). H
Then at a crossing $c_i$, edges $v_{i-1}$ and $v_{i}$ appear as lower hypotenuses of the ideal octahedron $o_i$ (see Figure \ref{fig:main_proofs} again). The Ptolemy coordinates on $o_i$ that we computed in Section~\ref{sec:main} show that up to sign
\[\dfrac{c'(v_{i-1})}{c'(v_i)}=\sigma(c_i).\]
It follows that
\begin{equation} 
	\sigma:=\prod_{i=1}^n \sigma(c_i)= \prod_{i=1}^n \dfrac{c'(v_{i-1})}{c'(v_i)} = \pm 1. 
\end{equation}
This proves the first claim of Theorem~\ref{thm:obs}.

\begin{claim}
\label{lem:pos} If $\sigma=1$, then $\Bcal : e(M)\rightarrow \psl$ lifts to a natural $(\sl,P)$-cocycle.
\end{claim}

\begin{proof} 		
	Recall Proposition \ref{prop:compatible} or \ref{prop:compatible_w} that the natural cocycle $\Bcal$ is induced from a collection of variables $p_i$ or $q_i \in \Cbb^\ast /\pmo$  satisfying certain relations. Even though  $p_i$ and $q_i$ are defined up to sign, the short-edge parameters are well-defined without sign-ambiguity, as they are even functions in $p_i$ and $q_i$. See the equations \eqref{eqn:z_short} and \eqref{eqn:w_short}. 
	In particular, we may lift $\Bcal(e)$ to $\widetilde{\Bcal}(e) \in \sl$ so that $\widetilde{\Bcal}(e) \in P$ for all short-edges $e$ of $M$. Note that then $\widetilde{\Bcal}$ satisfies the cocycle condition for all triangular faces of $M$.

	We now choose a lift $\widetilde{\Bcal}(e)$ of $\Bcal(e)$ for each  long-edge of $M$.  This is equivalent to choosing a lift $c: e(\Tcal)\rightarrow \Cbb^\ast$ of $c' : e(\Tcal)\rightarrow \Cbb^\ast/\pmo$.
	If a long-edge $e$ of $M$ joins two boundary spheres of $M$, then its long-edge parameter is an even function in $p_i$ and $q_i$. See Figures \ref{fig:z_ptolemy_pos}, \ref{fig:w_ptolemy_pos} and \ref{fig:w_ptolemy_neg}. Thus $c'(e)$ is in fact well-defined without sign-ambiguity. Letting $c(e):=c'(e)$, a lift $\widetilde{\Bcal}(e)$ of $\Bcal(e)$ is determined.

	We then consider a long-edge $v_i$ ($ 1 \leq i \leq n$). Recall that the corresponding arc $a_i$ of $D$ runs from the crossings $c_i$ to $c_{i+1}$. At each crossing in between $c_i$ and $c_{i+1}$, the long-edge $v_i$ appears twice as upper-hypotenuses. We enumerate the upper-hypotenuses at these intermediate crossings other than $v_i$ by $u_1,\ldots,u_m$ as in Figure \ref{fig:main_proofs}. For each $u_j$ $(1\leq j\leq m)$ there are two hexagonal faces, say $F_1$ and $F_2$, of $M$ containing both $v_i$ and $u_j$ (we have indicated $F_1$ and $F_2$ for $u_2$ in Figure \ref{fig:main_proofs}). The boundary of $F_1 \cup F_2$ consists of short-edges and long-edges joining two boundary spheres ($v_i$ and $u_j$ are common edges of  $F_1$ and $F_2$). In particular,  the $\widetilde{\Bcal}$-matrices for the boundary edges of $F_1 \cup F_2$ are already chosen and by construction $\widetilde{\Bcal}$ satisfies the cocycle condition for $F_1 \cup F_2$. This implies that whenever we choose a lift $\widetilde{\Bcal}(v_i)$, there is a unique lift $\widetilde{\Bcal}(u_j)$ so that $\widetilde{\Bcal}$ satisfies the cocycle condition for both $F_1$ and $F_2$.
	Applying the same argument for $u_m$ and $v_{i+1}$ (see the lower-hypotenuses at $c_{i+1}$), we deduce that a lift $\widetilde{\Bcal}(u_m)$ uniquely determines $\widetilde{\Bcal}(v_{i+1})$ so as to $\widetilde{\Bcal}$ satisfy the cocycle condition for each hexagonal face containing both $u_m$ and $v_{i+1}$.
	
	\begin{figure}[!h]
		\centering
		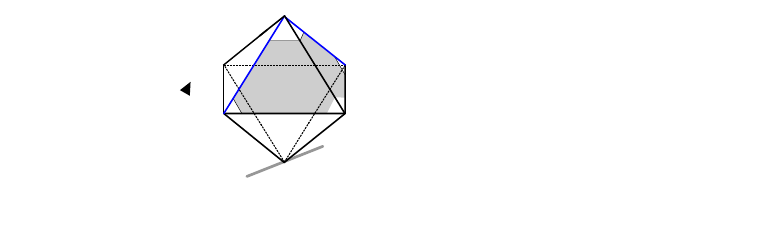
		\caption{An over-arc: faces $F_1$ and $F_2$ for $u_2$.}
		\label{fig:main_proofs}
	\end{figure}   	
	
	Regarding $c_1$ as an initial crossing, we choose any lift $\widetilde{\Bcal}(v_1)\in\sl$. The argument in the previous paragraph tells us that along the diagram, $\widetilde{\Bcal}(v_1)$ determines a lift $\widetilde{\Bcal}(e)$ uniquely  for all hypotenuses $e$ of the ideal octahedra $o_1,\ldots,o_n$ so that $\widetilde{\Bcal}$ satisfies the cocycle condition for all boundary faces of $\overline{o_1},\ldots\overline{o_n}$ except for lower hexagonal faces of $\overline{o_1}$.
	On the other hand, the uniqueness of such determination implies that    (see Figures \ref{fig:z_ptolemy_pos}, \ref{fig:w_ptolemy_pos},  \ref{fig:w_ptolemy_neg})
	\begin{equation} \label{eqn:ctilde}
		\dfrac{c(v_{i})}{c(v_{i+1})} = \sigma(c_{i+1}) \quad (1 \leq i\leq n-1).
	\end{equation}
	Moreover, the equation~(\ref{eqn:ctilde})  also holds for $i=n$, as $\sigma=1$. This implies that $\widetilde{\Bcal}$ satisfies the cocycle condition also for lower hexagonal faces of $\overline{o_1}$.
	
	We have chosen  $\widetilde{\Bcal}(e)$  for all boundary edges $e$ of $\overline{o_1},\ldots,\overline{o_n}$ to satisfy the cocycle condition for all boundary faces of $\overline{o_1},\ldots,\overline{o_n}$.
	Therefore, there is a lift $\widetilde{\Bcal}(e)$ for all internal edges $e$, edges  created to subdivide octahedra into ideal tetrahedra, so that $\widetilde{\Bcal} : e(M) \rightarrow \sl$ satisfies the cocycle condition for all faces of $M$. This proves that $\Bcal$ lifts to a natural $(\sl,P)$-cocycle $\widetilde{\Bcal}$.
\end{proof}

Claim \ref{lem:pos} shows that if $\sigma=1$, then $\rho$ lifts to an $(\sl,P)$-representation. In particular,  the obstruction class of $\rho$ is trivial. This proves the second part of Theorem~\ref{thm:obs} for $\sigma=1$.

\begin{claim}\label{lem:neg} If $\sigma= -1$, then $\rho$ admits an $\sl$-lifting $\widetilde{\rho}$ such that   $\textrm{tr}(\widetilde{\rho}(\lambda))=-2$ where $\lambda$ is a canonical longitude of $K$.
\end{claim} 
\begin{proof} We proceed the same construction of an $\sl$-lifting $\widetilde{\Bcal}$ of $\Bcal$ that we have done in the proof of Claim~\ref{lem:pos}. Then one can determine $\widetilde{\Bcal}(e)$ for all boundary edges $e$ of $\overline{o_1},\ldots,\overline{o_n}$ so that $\widetilde{\Bcal}$ satisfies the cocycle condition for all boundary faces of $\overline{o_1},\ldots,\overline{o_n}$, except lower hexagonal faces of $\overline{o_1}$. 
	In fact, as $\sigma=-1$, $\widetilde{\Bcal}$ fails the cocycle condition for two hexagonal faces $F_1$ and $F_2$ containing both $v_n$ and $u$, where $u$ is the lower-hypotenuse of $o_1$ other than $v_1$ and $v_{n}$ (see Figure \ref{fig:change_short}). 
	
	To make $\widetilde{\Bcal}$ be a cocycle, we change the sign of the matrix $\widetilde{\Bcal}(e_i)$ for the short-edge $e_i$ of $F_i$ $(i=1,2)$ near to the bottom vertex of $o_1$. If we are working on $\Tcal_4$, we also change the sign of $\widetilde{\Bcal}(e_3)$ for the short-edge $e_3$ described as in Figure \ref{fig:change_short}. Then one can easily check that $\widetilde{\Bcal}$ now satisfies the cocycle condition for all boundary faces of $\overline{o_1},\ldots,\overline{o_n}$.
	As in the proof of Claim~\ref{lem:pos}, we finally choose $\widetilde{\Bcal}(e)$ for internal edges  $e$ (created to subdivide octahedra into ideal tetrahedra) so that $\widetilde{\Bcal} :e(M) \rightarrow \sl$ satisfies the cocycle condition for all faces of $M$. We stress that $\widetilde{\Bcal}$ is a cocycle but not a natural $(\sl,P)$-cocycle, as $\widetilde{\Bcal}(e_i)$ is no longer in  $P$ for $1 \leq i \leq 3$.
	\begin{figure}[!h]
		\centering
		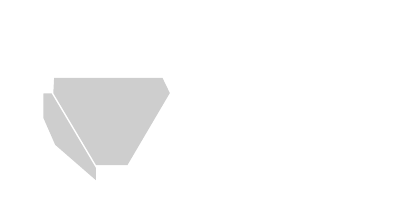
		\caption{Sign-change for the short-edges.}
		\label{fig:change_short}
	\end{figure}	
	
Let $\widetilde{\rho}:\pi_1(M)\rightarrow \sl$ be a representation induced from the cocycle $\widetilde{\Bcal}$. For a meridian $\mu$ of $K$ there is an edge-path of $M$ which is homotopic to $\mu$ and not passes the short-edges $e_1, e_2,$ and $e_3$. It follows that up to conjugation $\widetilde{\rho}(\mu) \in P$; in particular,  $\textrm{tr}\widetilde{\rho}(\mu)=2$. Similarly, for a (canonical) longitude $\lambda$ there is an edge-path of $M$ which is homotopic to $\lambda$ and passes one of $e_1, e_2, e_3$ exactly once. This proves that $\textrm{tr} \widetilde{\rho}(\lambda)=-2$.
\end{proof}

Combining Claim~\ref{lem:neg} with Proposition~2.2 in \cite{CYZ2019}, we obtain the second part of  Theorem \ref{thm:obs} for $  \sigma=-1$.

\subsection{Proof of Theorem~\ref{thm:longitude}}

For each crossing $c_i$ we denote by $d_i$ the diagonal of the bottom quadrilateral of $\overline{o_i}$ parallel to $D$. Precisely, the diagonal $d_i$ is homotopic to the composition of two short-edges $\cl{p} \ast \cl{m}$ in Figure~\ref{fig:z_short_edge}. 
From the equations \eqref{eqn:z_short} and \eqref{eqn:w_short} one can check that the sum of those two short-edge parameters is $\lambda(c_i)$.
This implies that if we let $\lambda_D$ be the blackboard framing longitude, then
\begin{equation*}
	\rho(\lambda_D) = \begin{pmatrix} 1 & \sum_{i=1}^n \lambda(c_i) \\ 0 & 1 \end{pmatrix}
\end{equation*}
Similarly, the meridian $\mu$ is homotopic to the composition of two short-edges $\cl{p} \ast \cl{o}$ and the sum of these short-edge parameters is $1$. That is,
\begin{equation*}
	\rho(\mu) = \begin{pmatrix} 1 &  1 \\ 0 & 1 \end{pmatrix}
\end{equation*}
Combining the above with the fact that $\lambda_D$ winds the knot $w(D)$ times, we deduce that the cusp shape of $\rho$ is $\sum \lambda(c_i) - w(D)$.

\begin{remark} One can use the diagonal of the top (instead of bottom) quadrilateral of $\overline{o}_i$ parallel to $D$. 
	Then a similar argument shows that the cusp shape is alternatively given by
	\[\displaystyle\sum_{i=1}^n \lambda'(c_i) -w(D)\]
	where $\lambda'(c_i)$ is given by
	\[\lambda'(c_i)=\dfrac{z_b-z_d}{z_c-z_a} \quad \textrm{or} \quad  \dfrac{w_aw_c- w_b w_d}{(w_a-w_b)(w_c-w_d)}\]
	for Figure \ref{fig:crossing_z}~(a) and~(b), or Figure \ref{fig:crossing_ww}~(a) and~(b).
\end{remark}

\section{Wirtinger generators of holonomy representation} \label{sec:wirtgen}
Based on the computation of a natural cocycle $\Bcal : e(M)\rightarrow \psl$ in Section~\ref{sec:main}, an explicit  Wirtinger generators of $\rho$ can be simply carried out: recall that for $\gamma \in \pi_1(M)$, $\rho(\gamma)$ is given by the product of $\Bcal$-matrices along an edge-path of $e(M)$ homotopic to $\gamma$ (see Proposition \ref{prop:hol_recv}). This section is devoted to present such an explicit algorithmic formula. 
We note that similar computations can be found in \cite[Section 5]{kim2018octahedral}, but the previoius formula was only valid under the condition that there should be no pinched octahedron, unlike the current formula for $\Tcal_5$ in this section.

For notational convenience,
we enumerate the crossings of a knot diagram $D$ by $c_i$ $(1\leq i \leq n)$ where the index of the crossings $c_i$ increases whenever we under-pass the crossing along the diagram as in Figure \ref{fig:main_proofs}.
We denote by $a_i$ the arc of $D$ running from $c_i$ to $c_{i+1}$  and denote by $\mu_i$ the Wirtinger generator winding the over-arc of a crossing $c_i$ respecting the orientation. See, for instance, Figure \ref{fig:trefoil_one}~(right). Let $e_i$ be the crossing-sign of $c_i$.

When segment variables around a crossing $c_i$ are given as in either Figure \ref{fig:crossing_z}~(a) or~(b), we define $M(c_i)\in\psl$ by
\begin{equation*}
	M(c_i):=
	\begin{pmatrix}
		\dfrac{z_a}{z_c} & \dfrac{z_b z_d(z_c-z_a)}{z_c^2(z_b-z_d)} \\[10pt]
		\dfrac{(z_a-z_c)(z_b-z_d)}{z_b z_d} & 2-\dfrac{z_a}{z_c} 
	\end{pmatrix}.
\end{equation*} 	
Similarly, when region variables around $c_i$ are given as in either Figure \ref{fig:crossing_ww}~(a) or~(b), we define $M(c_i) \in \psl$ by
\begin{equation*}
	M(c_i):=
	\begin{pmatrix}
		\dfrac{w_a-w_d}{w_b-w_c} & \dfrac{w_b w_d - w_a w_c}{(w_b-w_c)^2} \\[10pt]
		\dfrac{(w_a-w_b+w_c-w_d)^2}{w_a w_c -w_b w_d} & 2-\dfrac{w_a-w_d}{w_b-w_c} 
	\end{pmatrix}.
\end{equation*}
Then the $\rho$-image of the Wirtinger generators is inductively computed as follows.
\begin{theorem} \label{thm:w_wirtinger} For $1 \leq i \leq n$ we have 
	\begin{equation*}
		\rho(\mu_{i}^{e_i}) =	\rho(\mu_{1}^{e_{1}} \cdots \mu_{i-1}^{e_{i-1}})^{-1}
		\begin{pmatrix} 1 & \displaystyle\sum_{j=1}^{i} \lambda(c_j) \\ 0 &1 \end{pmatrix} M(c_i) \begin{pmatrix} 1 & \displaystyle\sum_{j=1}^{i} \lambda(c_j) \\ 0 &1 \end{pmatrix}^{-1} \rho(\mu_{1}^{e_{1}} \cdots \mu_{i-1}^{e_{i-1}})
	\end{equation*} 
\end{theorem}
\begin{proof} As in the proof of Theorem \ref{thm:longitude} we denote by $d_j$ the diagonal of the bottom quadrilateral of $\overline{o_j}$ parallel to the diagram $D$. Let $A_1$ be the composition $d_1\ast\cdots \ast d_i$, a path from $o_1$ to $o_i$ following $D$, and let $\underline{\mu_i}$ be a loop that (1) first follows the path $A_1$; (2) then winds the over-arc of the crossing $c_i$ respecting the orientation; (3) then returns through (the reverse of) $A_1$. See Figure \ref{fig:underline_mu}. 
	\begin{figure}[!h] 		
		\centering
		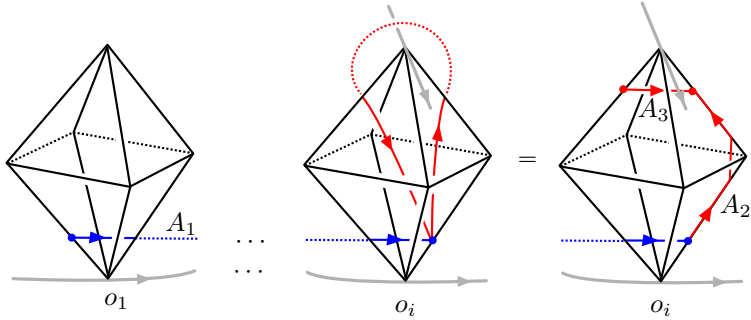
		\caption{The loop $\underline{\mu_{i}}$.}
		\label{fig:underline_mu}
	\end{figure}
	
	As we computed in the proof of Theorem \ref{thm:longitude}, the sum of short-edge parameters along the path $A_1$ is given by $\sum_{j=1}^{i} \lambda(c_j)$.
	For the second part (2) of the loop $\underline{\mu_i}$,
	using the parameters given in Section \ref{sec:five},
	the product of $\Bcal$-matrices along the edge-path $A_2\ast A_3 \ast A_2^{-1}$ depicted in Figure \ref{fig:underline_mu} is in fact the matrix $M(c_i)^{e_i}$. We thus obtain	\[ \rho(\underline{\mu_i}) = 	\begin{pmatrix} 1 & \displaystyle\sum_{j=1}^{i} \lambda(c_j) \\ 0 &1 \end{pmatrix} M(c_i)^{e_i} \begin{pmatrix} 1 & \displaystyle\sum_{j=1}^{i} \lambda(c_j) \\ 0 &1 \end{pmatrix}^{-1}.\]
	The theorem then follows from Lemma 5.1 in \cite{kim2018octahedral}, which says that the loop $\underline{\mu_i}$ is homotopic to $(\mu_1^{e_1} \cdots \mu_{i-1}^{e_{i-1}}) \mu_i (\mu_1^{e_1} \cdots \mu_{i-1}^{e_{i-1}})^{-1}$.
\end{proof}

\begin{example}\label{ex:trefoil_hol} We continue Example~\ref{ex:trefoil_obs} of the trefoil knot. Recall that we have $(e_1,e_2,e_3,e_4)=(1,1,-1,1)$ and $(\lambda(c_1),\lambda(c_2),\lambda(c_3),\lambda(c_4))=(-\frac{9}{2},\frac{4}{3},-1,\frac{1}{6})$. A simple computation gives
	\[ M(c_1) = \begin{pmatrix} -2 & 9 \\ -1 & 4 \end{pmatrix},\  M(c_2) = \begin{pmatrix} \frac{1}{3} & \frac{4}{9} \\[3pt] -1 & \frac{5}{3} \end{pmatrix},\ M(c_3) = \begin{pmatrix} 1 & -1 \\ 0 & 1 \end{pmatrix},\ M(c_4) = \begin{pmatrix} \frac{3}{2} & \frac{1}{4} \\[3pt] -1 & \frac{1}{2} \end{pmatrix}.\]
	Applying Theorem \ref{thm:w_wirtinger}, the $\rho$-image of the Wirtinger generators is computed as follow (see Figure \ref{fig:trefoil_one}~(right)):
	
	\begin{align*}
		\rho(\mu_1^{e_1}) & = 	\arraycolsep=4pt\begin{pmatrix} 1  & -\frac{9}{2} \\ 0 &1 \end{pmatrix} M(c_1) \begin{pmatrix} 1 & -\frac{9}{2} \\ 0 &1 \end{pmatrix}^{-1}=
		\begin{pmatrix} \frac{5}{2} &  \frac{9}{4} \\[3pt]-1 & -\frac{1}{2} 		
		\end{pmatrix} , \\
	\rho(\mu_2^{e_2}) &= 	\arraycolsep=4pt\rho(\mu_1^{e_1})^{-1} 	\begin{pmatrix} 1  & -\frac{19}{6} \\ 0 &1 \end{pmatrix}M(c_2)	\begin{pmatrix} 1  &-\frac{19}{6} \\ 0 &1 \end{pmatrix}^{-1}\rho(\mu_1^{e_1})=\begin{pmatrix} 1 &  1 \\[3pt] 0 & 1
	\end{pmatrix} , \\
\rho(\mu_3^{e_3})&=	\arraycolsep=4pt\rho(\mu_1^{e_1}\mu_2^{e_2})^{-1}\begin{pmatrix} 1  & -\frac{25}{6} \\ 0 &1 \end{pmatrix}M(c_3)	\begin{pmatrix} 1  &-\frac{25}{6} \\ 0 &1 \end{pmatrix}^{-1}\rho(\mu_1^{e_1}\mu_2^{e_2})=\begin{pmatrix} -\frac{1}{2} &  -\frac{9}{4} \\[3pt] 1 & \frac{5}{2}
\end{pmatrix},\\
\rho(\mu_4^{e_4})&=	\arraycolsep=4pt\rho(\mu_1^{e_1}\mu_2^{e_2}\mu_3^{e_3})^{-1}\begin{pmatrix} 1  & -4 \\ 0 &1 \end{pmatrix}M(c_4)	\begin{pmatrix} 1  &-4 \\ 0 &1 \end{pmatrix}^{-1}\rho(\mu_1^{e_1}\mu_2^{e_2}\mu_3^{e_3})=\begin{pmatrix} \frac{3}{2} &  \frac{1}{4} \\[3pt] -1 & \frac{1}{2}
\end{pmatrix} .
	\end{align*}
\end{example}

\section{Acknowledgment}
HK is	Supported by Basic Science Research Program through the NRF of Korea funded by the Ministry of Education (NRF-2018R1A2B6005691). SK is supported  by Basic Science Research Program through the National Research Foundation of Korea(NRF) funded by the Ministry of Education(2022R1I1A1A01063774) and
by  the Institute for Basic Science (IBS-R003-D1).

\bibliographystyle{alpha}
\bibliography{biblog}
\end{document}